\documentclass{article}
\usepackage[utf8]{inputenc}

\usepackage{fullpage}
\usepackage{amsmath}
\usepackage{amssymb}
\usepackage{graphicx}
\usepackage{subcaption}
\usepackage{comment}
\usepackage{soul, xcolor}
\usepackage{authblk}
\usepackage{algorithm}
\usepackage{algorithmic}
\usepackage{makecell} 

\usepackage[colorlinks]{hyperref} 

\usepackage[
    backend=biber,
    style=numeric,
    sortcites=true,
    sorting=none,
    hyperref=true,
    backref,
    natbib,
]{biblatex}

\hypersetup{
    citecolor = blue,
}

\graphicspath{{Figures/}}

\title{A data driven heuristic for rapid convergence of Scheduled Relaxation Jacobi schemes}

\author[,a]{Mohammad Shafaet Islam\thanks{Corresponding author \vfill \ \ \textit{Email addresses}: \texttt{moislam@mit.edu} (Mohammad Shafaet Islam), \texttt{qiqi@mit.edu} (Qiqi Wang)}}
\author[a]{Qiqi Wang}

\affil[a]{\footnotesize{\textit{Department of Aeronautics and Astronautics, Massachusetts Institute of Technology, 77 Massachusetts Avenue, Cambridge, MA 02139, USA}}}

\date{}

\begin{document}

\maketitle

\begin{center}
\rule{\textwidth}{.5pt}  
\end{center}
\begin{abstract}
The Scheduled Relaxation Jacobi (SRJ) method is a viable candidate as a high performance linear solver for elliptic partial differential equations (PDEs). The method greatly improves the convergence of the standard Jacobi iteration by applying a sequence of $M$ well-chosen overrelaxation and underrelaxation factors in each cycle of the algorithm to effectively attenuate the solution error. In previous work, optimal SRJ schemes (sets of relaxation factors) have been derived to accelerate convergence for specific discretizations of elliptic PDEs. In this work, we develop a family of SRJ schemes which can be applied to solve elliptic PDEs regardless of the specific discretization employed. To achieve favorable convergence, we train an algorithm to select which scheme in this family to apply at each cycle of the linear solve process, based on convergence data collected from applying these schemes to the one-dimensional Poisson equation. The automatic selection heuristic that is developed based on this limited data is found to provide good convergence for a wide range of problems. 


\end{abstract}
\begin{center}
\rule{\textwidth}{.5pt}  
\end{center}

\section{Introduction}

The solution of large linear systems of equations is an important problem in computational science and engineering, and is of great practical interest to scientists and engineers who are modeling physical phenomena such as fluid flow \cite{cfd2030} or electromagnetics \cite{electromagnetics-book}. These phenomena are usually governed by partial differential equations (PDEs) which may not be amenable to analytic solution but can usually be solved numerically. The numerical solution of PDEs by a numerical method leads to a large sparse linear system of equations which must be solved efficiently.

A host of methods have been developed in order to solve these linear systems, usually classified into direct and iterative methods \cite{Bau97}. Direct methods solve linear systems exactly but may become intractable as the number of degrees of freedom grows. Iterative methods have become popular in recent years as computing capability has continued to improve due to advances in hardware. Stationary iterative methods such as Jacobi iteration \cite{jacobi} and Gauss-Seidel were the first iterative methods used for solving large linear systems, and utilize relaxation steps to remove components of the residual vector in order to converge towards the exact solution \cite{Saad03}. However, the number of iterations required for convergence grows rapidly with the size of the system \cite{ParallelIterativeAlgorithms2007}. Krylov subspace methods such as conjugate gradient \cite{CG} and GMRES \cite{gmres} achieve faster convergence compared to stationary iterative methods and are therefore the most prominent iterative methods for solving sparse linear systems of equations. Despite this, Jacobi iteration is exceptionally well suited to implementation on the latest high performance computers and can achieve good performance due to its highly parallel nature and remarkable simplicity. Furthermore, Jacobi iteration is commonly applied as a smoother in multigrid solvers which exhibit fast convergence \cite{Briggs2000}. Therefore, theoretical improvements to the convergence of Jacobi iteration can augment its usability on high performance systems, and potentially make it a more viable method for solving large linear systems compared to the Krylov subspace methods, whose performance may be inhibited on large hierarchical high performance computing systems due to the reduction operations (e.g. dot products) that are required at every step of these algorithms.

There has been recent work to improve the convergence of Jacobi iteration. The Scheduled Relaxation Jacobi (SRJ) method developed by Yang and Mittal accelerates the convergence of the standard Jacobi iteration by applying prescribed relaxation factors in a cyclic fashion \cite{Yang14}. In each cycle of the algorithm, a fixed number of relaxed Jacobi iterations $M$ are performed with $P$ distinct predetermined factors which depend on the SRJ scheme one chooses to use. Each scheme is characterized by the number of distinct relaxation values used in the cycle (denoted by $P$), the specific relaxation factors used (listed in the vector $\vec{\Omega} = [\omega_1, \omega_2, ..., \omega_P]$ where $\omega_1 > \omega_2 > ... > \omega_P$), and the frequency with which each factor $\omega_i$ is applied (denoted by $q_i$ and stored in the vector $\vec{Q} = [q_1, q_2, ..., q_P]$). For a given $P$, one can derive a scheme characterized by the distinct relaxation factors $\vec{\Omega}$ and their frequency of use $\vec{Q}$ which results in optimal convergence. This is done by solving a min-max optimization problem in which the maximum possible amplification factor associated with the scheme is minimized. Yang and Mittal derive a variety of schemes for $P = 2,3,4,5$ and derive parameters for different grid sizes $N$. Typical schemes involve a few overrelaxation steps followed by many underrelaxation steps so that the overall effect of the iterations attenuates certain modes. For a certain $P = 7$ scheme, they observe a 190 times speedup in convergence for the 2D Laplace and Poisson equations relative to standard Jacobi iteration. While schemes corresponding to larger $P$ would likely provide even faster convergence,  deriving such schemes is difficult as the system of equations for the scheme parameters grows stiffer as $P$ and $N$ are increased. Adsuara et al. proposed a number of algebraic simplifications in order to make it easier to derive more complicated schemes for $P$ up to 15 and $N = 2^{15}$. They report even higher speedups than Yang and Mittal, nearly to a factor of 1000 for large $P$ and problem size $N$ \cite{Adsuara1}. In \cite{Adsuara2}, Adsuara et al. propose SRJ schemes in which each relaxation factor is only used once per cycle. They reason that although a scheme with sufficiently large $P$ may greatly reduce the solution residual, the overall number of iterations associated with the cycle $M$ may be extremely large if each distinct relaxation factor is used many times. In this case, the SRJ scheme may not necessarily outperform Jacobi iteration. They introduce the Chebyshev Jacobi method (CJM) and find that Yang's original min-max optimization problem for deriving an optimal SRJ scheme becomes a simpler problem of solving for the roots of a scaled Chebyshev polynomial, the reciprocal of which are the desired relaxation parameters. A CJM scheme can be derived for a given problem size based on the minimum and maximum wavenumbers. For a fixed number of iterations, their schemes with distinct relaxation factors converge faster relative to the original schemes by Yang and Mittal where relaxation parameters are repeatedly used. Specifically, SRJ schemes where $P = M$ are most effective for fast convergence.

Both the SRJ schemes developed in \cite{Yang14} and \cite{Adsuara1}, and the CJM schemes developed in \cite{Adsuara2} are specific to a problem size $N$. This presents the advantage that the schemes are tailored to the specific problem of interest and will provide optimal convergence. However, this also presents certain difficulties in terms of practical implementation. For example, if a practitioner is interested in solving a variety of problems or even a single problem under several discretizations, they will be required to construct a new scheme for each case, which can be a time-consuming process. Furthermore, given a fixed discretization, determining the value of $M$ which provides the optimal convergence can also be tedious especially when this must be done for many problems. This work aims to address both of these practical aspects by introducing a more general approach. First, we consider a family of SRJ schemes that can be applied to all symmetric linear systems which would converge via Jacobi iteration. Within our family of schemes, determining the best scheme for a given problem may be difficult. In some cases, it may be possible to determine an asymptotic convergence rate associated with different SRJ schemes for solving a linear system, so that one can choose a scheme with a relatively large convergence rate. However, for many large scale problems computing these convergence rates may not be feasible, in which case a good scheme must be found through brute force experimentation for each linear system under consideration. As an alternative approach, we develop a data driven heuristic to determine which scheme to use at a given step of the solution process so that we can obtain rapid convergence for a variety of problems. This heuristic is trained on convergence data which we collect from applying the SRJ schemes to select training problems. The approach avoids the need to develop new schemes for different discretizations of an elliptic PDE, and also avoids the need to determine an optimal $M$ as scheme selection from our family of schemes will be determined according to this data informed heuristic. 

The remainder of this paper is organized as follows. Section 2 illustrates an approach for developing our family of SRJ schemes. These schemes are not restricted to specific problem sizes as those in \cite{Yang14} and \cite{Adsuara1} so they may be used to solve a variety of symmetric linear systems such as those arising from discretization of elliptic PDEs (the only requirement is that the original systems could be solved by the standard Jacobi iterative method). Section 3 presents the data driven approach used to inform the automatic selection process for selecting which scheme to use in a given SRJ cycle. A simple heuristic is developed to select a scheme for the next cycle. Section 4 shows the performance of the SRJ schemes with this data based automatic selection heuristic for both in-sample matrices from which we collected data, as well as many out-of-sample test matrices which are progressively more and more distinct from the in-sample matrices. These sample linear systems arise from discretizations of elliptic PDEs on structured domains, as well as unstructured meshes which have not been explored before in the context of the Scheduled Relaxation Jacobi method. Section 5 provides concluding remarks and a vision for SRJ on high performance computers. Our hope is that practitioners can utilize the SRJ schemes developed here (based on the theory presented in Section 2), and apply our scheme selection heuristic to solve their problems of interest without the need to tune the schemes or any additional parameters. 

\section{Derivation of SRJ schemes}

We present an approach to derive a set of relaxation factors (an SRJ scheme) which will improve convergence of Jacobi iteration when applied to a linear system of equations $Ax = b$, $A \in \mathbb{R}^{n \times n}$, $x \in \mathbb{R}^{n}$, $b \in \mathbb{R}^{n}$. We begin by defining the original Jacobi iterative method. This involves a matrix splitting of $A$ into the following form
\begin{equation}
    A = D + L + U 
\end{equation}
where $D$ is a diagonal matrix containing the diagonal entries of $A$, and $L$ and $U$ are matrices containing the lower and upper portions of $A$ below and above the diagonal. The linear system can be written as the following fixed point update from step $n$ to step $n+1$
\begin{equation}
    x^{(n+1)} = \underbrace{-D^{-1} (L + U)}_{B_\text{J}} x^{(n)} + D^{-1} b
    \label{eqn:jacobi}
\end{equation}
Equation \eqref{eqn:jacobi} represents the Jacobi iterative update for solving a linear system of equations where $B_\text{J} = -D^{-1} (L+U) $ is the iteration matrix associated with the Jacobi update. The rate of convergence of the method depends on the spectral radius of the iteration matrix $B_\text{J}$, which must be less than 1 for Jacobi iteration to converge \cite{Golub2007}. We now consider the weighted Jacobi iteration which is given by the following update equation
\begin{equation}
    x^{(n+1)} = \omega \left[ B_{\text{J}} x^{(n)} + D^{-1} b \right] + (1 - \omega) x^{(n)}
    \label{eqn:relaxed-jacobi}
\end{equation}
and involves weighting the original update equation \eqref{eqn:jacobi} by some relaxation factor $\omega$. Setting $\omega < 1$ is known as underrelaxation while setting $\omega > 1$ is referred to as overrelaxation (applying overrelaxation alone is known to cause Jacobi iteration to diverge). The update equation \eqref{eqn:relaxed-jacobi} can also be written as
\begin{equation}
    x^{(n+1)} = \left[ (1-\omega) I + \omega  B_{\text{J}} \right] x^{(n)} + \omega D^{-1} b 
    \label{eqn:relaxed-jacobi-simplified}
\end{equation}
We define the error vector at step $n$ as $e^{(n)} \equiv x^{(n)} - x$ where $x$ is the exact solution to the linear system. The exact solution satisfies the update equation exactly as follows
\begin{equation}
    x = \left[ (1-\omega) I + \omega  B_{\text{J}} \right] x + \omega D^{-1} b 
    \label{eqn:exact-update}
\end{equation}
Subtracting Equation \eqref{eqn:exact-update} from \eqref{eqn:relaxed-jacobi-simplified} yields an equation for the evolution of the error vector from one step to the next as given by Equation \eqref{eqn:error-evolution-equation}.
\begin{equation}
    e^{(n+1)} = \underbrace{\left[ (1-\omega) I + \omega  B_{\text{J}} \right]}_{B_{\omega}} e^{(n)} 
    \label{eqn:error-evolution-equation}
\end{equation}
The accumulation of error is based on the amplification matrix $B_{\omega}$. Convergence of weighted Jacobi is guaranteed if the matrix $B_{\omega}$ has a spectral radius less than 1.

We now consider an iteration scheme where $M$ iterations of the relaxed Jacobi method with distinct $\omega$ are performed. Let these iterations comprise one cycle of the SRJ method, and denote the overall amplification matrix associated with these $M$ iterations by  $B_{\text{SRJ}}$. Also denote the amplification matrix associated with each individual iteration by $B_{\omega_{i}}$. If $e^{(n)}$ and $e^{(n+1)}$ represent the error prior to and after a cycle of $M$ iterations, then the error accumulates as follows
\begin{equation}
    e^{(n+1)} = B_{\text{SRJ}} e^{(n)}  = \prod_{i=1}^{M} B_{\omega_i} e^{(n)}  = \prod_{i=1}^{M} \left[ (1-\omega_i) I + \omega_i  B_{\text{J}} \right] e^{(n)} 
\end{equation}
The amplification of the error at each cycle is related to the eigenvalues of the matrix $B_{\text{SRJ}}$. In order for the error to decay from one SRJ cycle to the next, the spectral radius of the SRJ iteration matrix must be less than 1. In other words, all eigenvalues of $B_{\text{SRJ}}$ must have magnitude less than 1. It can be shown that the eigenvalues of $B_{\text{SRJ}}$ (which we denote by $\lambda_{\text{SRJ}}$) are related to the eigenvalues of $B_{\text{J}}$ (which we denote by $\lambda_{\text{J}}$). If $v_{j}$ denotes an eigenvector of both $B_{\text{SRJ}}$ and $B_{\text{J}}$, it follows that
\begin{align}
    B_{\text{SRJ}} v_{j} &= \prod_{i = 1}^{M} \left[ (1-\omega_{i}) I + \omega_{i} B_{\text{J}} \right] v_{j} \\ 
    &= \prod_{i = 1}^{M} \left[ (1-\omega_{i}) v_{j} + \omega_{i} \lambda_{\text{J}} v_{j} \right]  \\ 
    &= \prod_{i = 1}^{M} \left[ (1-\omega_{i}) + \omega_{i} \lambda_{\text{J}} \right] v_{j} 
    \label{eqn:eigen-srj-1}
\end{align}
Additionally, by definition it is true that
\begin{equation}
     B_{\text{SRJ}} v_{j} = \lambda_{\text{SRJ}} v_{j}
     \label{eqn:eigen-srj-2}
\end{equation}
Given Equations \eqref{eqn:eigen-srj-1} and \eqref{eqn:eigen-srj-2} are true, we obtain the following relationship between the eigenvalues of the SRJ iteration matrix and those of the Jacobi iteration matrix
\begin{equation}
    \lambda_{\text{SRJ}} = G_M(\lambda_\text{J})\;,\quad\mbox{where}\quad G_M(\lambda):= \prod_{i=1}^{M} \left[(1-\omega_i) + \omega_i \lambda\right]
    \label{eqn:amplification-M-scheme}
\end{equation}
According to Equation \eqref{eqn:amplification-M-scheme}, the eigenvalues of $B_{\text{SRJ}}$ are found by evaluating an $M$-degree polynomial at the eigenvalues of $B_\text{J}$, which need to be in $(-1,1)$
for Jacobi iteration to converge. We call this $M$-degree polynomial the amplification polynomial associated with the SRJ scheme, and denote it by $G_M$. In order to design an SRJ scheme with $M$ relaxation parameters, we only need to design a $M$-order polynomial $G_M$ and ensure it has $M$ real roots.  The value of $G_M(\lambda)$ should lie in $(-1,1)$ when $\lambda \in (-1,1)$, in order for $B_{\text{SRJ}}$ to have a spectral radius less than 1 when the Jacobi iteration matrix $B_\text{J}$ has a spectral radius less than 1. In fact, one can deliberately construct the amplification polynomial so that it converges more rapidly than Jacobi for ranges of $\lambda_{\text{J}}$. 

\begin{figure}[htbp!]
    \centering
    \includegraphics[width=0.5\textwidth]{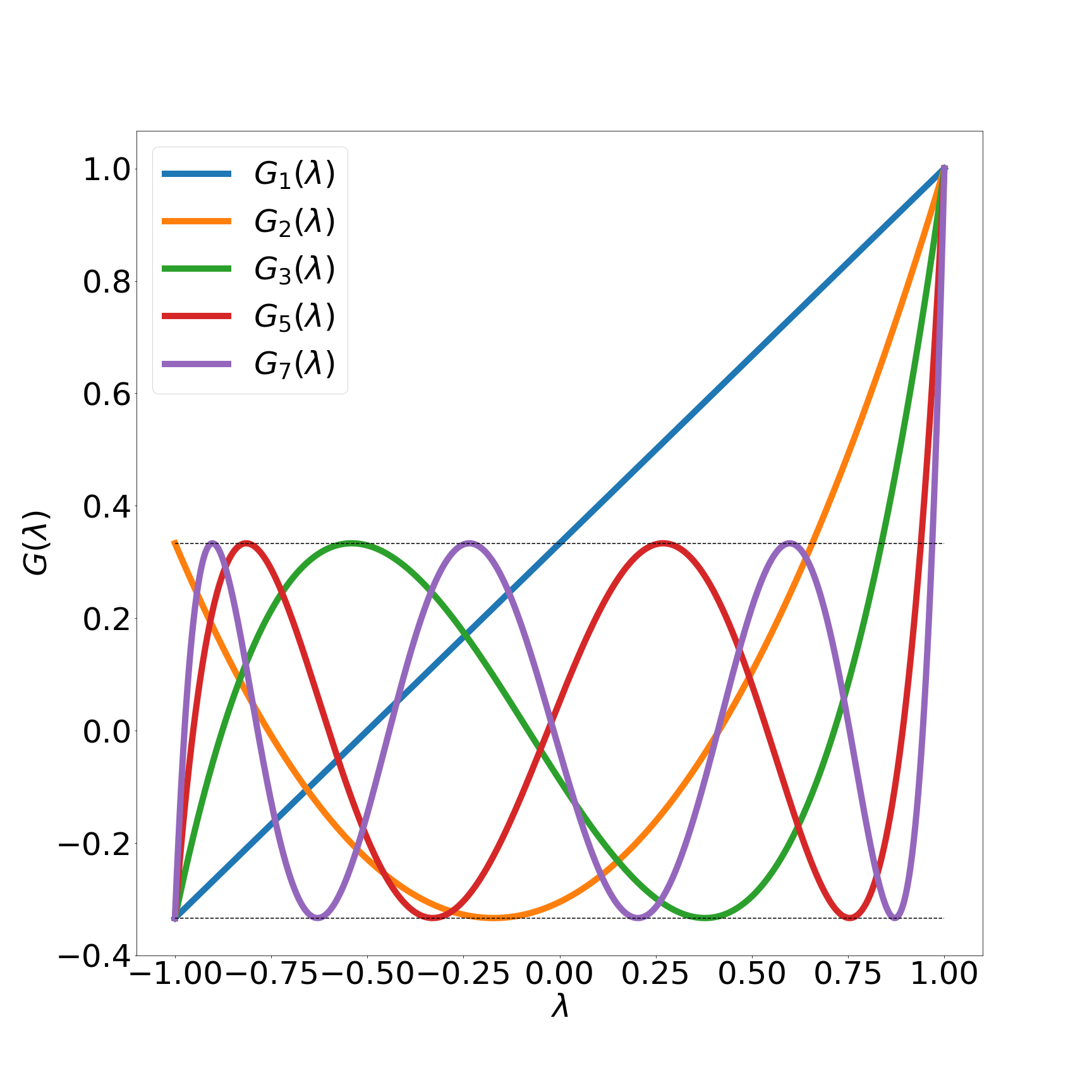}
    \caption{Amplification polynomials $G_{M}(\lambda)$ for $M = 1,2,3,5,7$. The polynomials are bounded by $\frac{1}{3}$ for some region within $\lambda \in (-1,1)$ which grows as $M$ increases. Applying the SRJ schemes corresponding to these amplification polynomials results in many of the solution error eigenmodes being attenuated by at least a factor of $\frac{1}{3}$.}
    \label{fig:amplification}
\end{figure}

We construct a sequence of amplification polynomials of increasing order $M$.  The resulting SRJ scheme has convergence properties suitable for matrices of different stiffness. Our amplification polynomials are shown in Figure \ref{fig:amplification} for $M = 1,2,3,5,7$. These polynomials are constructed such that they are bounded for the widest range possible in $(-1,1)$ by some bounding value. For $M = 1$, the polynomial is a straight line and attenuates all modes associated with eigenvalues between $(-1,0)$ by at least a factor of $\frac{1}{3}$. This bounding value corresponds to a Jacobi relaxation scheme with relaxation factor $\frac{2}{3}$, which is a popular choice as a smoother for multigrid methods because it decays the higher order modes of the error by a factor of $\frac{1}{3}$ each iteration 
\cite{Briggs2000}. As the degree of the polynomial $M$ grows, the range of eigenvalues over which the amplification is bounded also grows. Within this range, the amplification polynomial is actually a scaled version of a Chebyshev polynomial. These polynomials have the property that all extrema have the same absolute value, so they are useful as polynomials which are bounded in some interval. We denote the $M$th Chebyshev polynomial by $T_{M}$.

Our amplification polynomials are scaled versions of the Chebyshev polynomials. One can derive the relationship between the amplification polynomials and the Chebyshev polynomials (done in Appendix A), which is 
\begin{equation}
    G_{M}(\lambda) = \frac{T_{M} (f(\lambda)) }{3} \;,\quad\mbox{where}\quad f(\lambda):= \frac{(\lambda^* + 1)\lambda + (\lambda^*-1)}{2}
    \label{eqn:transformation}
\end{equation}
where $\lambda^*$ satisfies $T_M (\lambda^*) = 3$. Table \ref{tab:amplification-polynomials} shows the amplification polynomials for $M = 1,2,3,5$ (plotted in Figure \ref{fig:amplification}) as well as the corresponding Chebyshev polynomials. The maximum eigenvalue $\lambda_{\text{max}}$ for which the amplification polynomial is bounded by $\frac{1}{3}$ is also shown (and grows closer to 1 with polynomial order $M$). The amplification polynomial is always bounded at the minimum eigenvalue of $\lambda_{\text{min}} = -1$ for all $M$.

\begin{table}[htbp!]
    \caption{The amplification polynomials $G_{M}$ and the maximum $\lambda$ for which the polynomials are bounded by $\frac{1}{3}$, for $M = 1,2,3,5$. The corresponding Chebyshev polynomials $T_{M}(\lambda)$ are also shown.}
    \centering
    \begin{tabular}{|c|l|l|c|}
        \hline
        $M$ & $T_{M}(\lambda)$ & $G_{M}(\lambda)$ & $\lambda_{\text{max}}$  \\
        \hline
        1 & $\lambda$ & $\frac{2}{3}\lambda + \frac{1}{3}$ & 0.0 \\
        \hline
        2 & $2\lambda^2-1$ & $\frac{2}{3} (1.2071 \lambda + 0.2071)^2 - \frac{1}{3}$ & 0.6569 \\
        \hline
        3 & $4\lambda^3-3\lambda$ & $\frac{4}{3} (1.0888 \lambda + 0.0888)^3 - (1.0888 \lambda + 0.0888)$ & 0.8368 \\
        \hline
        5 & $16\lambda^5 - 20\lambda^3 + 5\lambda$ & $\frac{16}{3}(1.0314 \lambda + 0.0314)^5 - \frac{20}{3}(1.0314 \lambda + 0.0314)^3 + \frac{5}{3}(1.0314 \lambda + 0.0314)$ & 0.9391 \\
        \hline
    \end{tabular}
    \label{tab:amplification-polynomials}
\end{table}

The SRJ scheme with $M$ relaxation factors can be derived given the $M$-degree amplification polynomial. Given the $j$th root of $G_{M}(\lambda)$ which we denote by $\lambda^{r}_j$, it is true that
\begin{equation}
    G_{M} \left( \lambda^{r}_j \right) = 0 \rightarrow \prod_{i=1}^{M} \left[ (1-\omega_i)+ \omega_i  \lambda^{r}_{j} \right] = 0 \rightarrow 
    (1-\omega_j)+ \omega_j \lambda^{r}_j = 0 
\end{equation}
The $j$th relaxation factor in the SRJ scheme is related to the $j$th root of the polynomial by Equation \eqref{eq:optimal-relaxed}
\begin{equation}
    \omega_j = \frac{1}{1-\lambda^{r}_{j}} \ , \ \forall \ j \in [1,M]
    \label{eq:optimal-relaxed}
\end{equation}

We tabulate the relaxation factors associated with the SRJ schemes for size $M = 1,2,3,5,7$ in Table \ref{tab:srj-schemes}. Applying the $M$ relaxation factors given in Equation \eqref{eq:optimal-relaxed} in an SRJ cycle results in the solution error eigenmodes being amplified by $G_{M}(\lambda)$. Each scheme results in a different amplification of the eigenmodes. Larger schemes attenuate modes corresponding a larger range of eigenvalues but are more computationally intensive to execute as they involve more iterations. 

\begin{table}[htbp!]
    \caption{Relaxation factors associated with SRJ schemes derived for $M = 1,2,3,5,7$}
    \centering
    \begin{tabular}{|c|l|}
        \hline
        $M$ & SRJ scheme parameters \\
        \hline
        1 & 0.66666667 \\
        \hline
        2 & 1.70710678, 0.56903559 \\
        \hline
        3 & 3.49402108, 0.53277784, 0.92457411 \\
        \hline
        5 & 9.23070105, 0.51215173, 0.97045899, 0.62486988, 2.1713295 \\
        \hline
        7 & 17.84007924, 0.50624677, 0.9845549, 1.69891732, 0.56014439, 4.06304526, 0.69311375\\
        \hline
    \end{tabular}
    \label{tab:srj-schemes}
\end{table}

We can define a procedure to compute the $M$ order SRJ scheme for arbitrary $M$. The SRJ scheme is related to the roots of the $M$ order amplification polynomial, which is a linearly transformed Chebyshev polynomial. Therefore, given the roots of the Chebyshev polynomial, a transformation can be applied to obtain the roots of the amplification polynomial (which can then be used to obtain the SRJ scheme relaxation parameters) . It is true that the amplification and Chebyshev polynomials are related by $G_{M}(\lambda) = \frac{T_M(f(\lambda))}{3}$ where $f(\lambda)$ is the transformation given in Equation \eqref{eqn:transformation}, so it follows that $G_{M}(f^{-1}(\lambda)) = \frac{T_M(\lambda)}{3}$. Denote the transformation $f^{-1}$ as $g$. Given the roots of $T_M$ denoted by $x^{r}_{j}$, the roots of $G_{M}$ are $\lambda^{r}_j$ where
 \begin{equation}
    \lambda^{r}_{j} = g(x^{r}_{j}) \;,\quad\mbox{where}\quad g(x):=  \frac{2}{\lambda^*+ 1} x + \frac{1 - \lambda^*}{1+\lambda^*}
    \label{eqn:transformation-inverse}
\end{equation}
As before, $\lambda^*$ satisfies $T_M(\lambda^*) = 3$. This gives a relationship between the roots of the Chebyshev polynomial and those of the amplification polynomial. In summary, determining the $M$ order SRJ scheme involves the following three steps:
\begin{enumerate}
    \item Find the roots of the $M$-degree Chebyshev polynomial $T_{M}$ denoted by $x^{r}_{j}$.
    \item Given the roots $x^{r}_{j}$ of the Chebyshev polynomial $T_M$, the roots of the amplification polynomial $G_{M}$ are $\lambda^{r}_{j} = g \left(x^{r}_{j} \right)$, where $g$ is the transformation defined in Equation \eqref{eqn:transformation-inverse}.
    \item Given the roots of the amplification polynomial $\lambda^{r}_{j}$, solve for the corresponding relaxation factors $\omega_j$ using Equation \eqref{eq:optimal-relaxed}.
\end{enumerate}
The three steps above can be combined into a single expression for the SRJ scheme of length $M$. Given the $M$-degree Chebyshev polynomial $T_{M}$, its roots $x^{r}_{j}$, and the argument $\lambda^{*}$ which satisfies $T_{M}(\lambda^{*}) = 3$, the relaxation factors of the SRJ scheme of length $M$ are given by
\begin{equation}
    \omega_{j} = \frac{\lambda^{*}+1}{2 \left(\lambda^{*}-x^{r}_{j} \right)}, \ j \in \left[0,M-1 \right]
    \label{eqn:srj-relaxation-parameters}
\end{equation}

Equation \eqref{eqn:srj-relaxation-parameters} gives the SRJ scheme relaxation factors which can be employed in the relaxed Jacobi iteration (i.e. Equation \eqref{eqn:relaxed-jacobi-simplified}) to improve convergence when solving linear systems arising from discretization of elliptic PDEs. A family of SRJ schemes can be derived, with each scheme being associated with a different length $M$. The SRJ schemes developed here can be expressed in the context of other polynomial acceleration methods which have previously been developed \cite{varga1962iterative}. For example, the SRJ schemes have the same effect as applying the Chebyshev semi-iterative method \cite{golub1961chebyshev} to accelerate Jacobi iteration, if the lower and upper eigenvalue bounds are specified as -1 and $\lambda_{\text{max}}$ respectively for a given $M$. The relaxed iterations in Equation \eqref{eqn:relaxed-jacobi-simplified} associated with the SRJ method can also be expressed in terms of Richardson iterations \cite{richardson}. In this case, the relaxation factors used in Richardson are the SRJ relaxation factors scaled by the matrix $D^{-1}$ (as explained in \cite{Adsuara2}). 

There are two useful metrics for comparing the performance of SRJ schemes of differing length $M$ for solving a linear system. One metric is the asymptotic convergence rate per iteration associated with the SRJ scheme for solving a linear system of interest. This convergence rate can be computed as follows
\begin{equation}
    \text{Asymptotic Convergence Rate/Iteration} = -\frac{\log \rho \left(B_{\text{SRJ}} \right)}{M} 
    \label{eqn:convergence-rate}
\end{equation}
and depends on the spectral radius of the SRJ iteration matrix (denoted by $\rho \left(B_{\text{SRJ}}\right)$) which can be computed based on Equation \eqref{eqn:spectral-radius-srj} as follows
\begin{equation}
    \rho \left(B_{\text{SRJ}} \right) = \max |\lambda_{\text{SRJ}}| = \max |G_{M} \left(\lambda_{\text{J}} \right)| 
    \label{eqn:spectral-radius-srj}
\end{equation}
According to Equation \eqref{eqn:convergence-rate}, one can obtain the asymptotic convergence rate per iteration of the SRJ scheme based on the spectral radius of the SRJ iteration matrix. Obtaining this spectral radius requires computing the SRJ iteration matrix eigenvalues, which can be found according to Equation \eqref{eqn:spectral-radius-srj} by evaluating the amplification polynomial $G_{M}$ associated with the SRJ scheme at the Jacobi iteration matrix eigenvalues corresponding to our linear system of interest. For a given problem, the scheme which yields the largest per iteration asymptotic convergence rate is expected to provide the fastest convergence. Another metric for comparing SRJ schemes is their ability to solve stiff problems. The ability of a SRJ scheme to converge for a stiff problem can be characterized by the slope of its corresponding amplification polynomial at $\lambda = 1$ (i.e. $G_{M}^{'}(\lambda = 1)$). If the slope is larger, the spectral radius of the SRJ scheme for the stiff problem is likely to deviate further away from 1, corresponding to a faster asymptotic convergence rate. Schemes associated with larger $M$ are better able to handle stiff problems. This is because the amplification polynomials associated with larger $M$ bound a larger portion of the $\lambda \in (-1,1)$ region. A stiff problem would have Jacobi iteration matrix eigenvalues close to 1 (so the spectral radius of the Jacobi iteration matrix is close to 1), so schemes associated with larger $M$ are better able to "capture" these eigenvalues and bound their associated amplification (so the spectral radius of the associated SRJ iteration matrix is further away from 1). 

The SRJ schemes based on Equation \eqref{eqn:srj-relaxation-parameters} can be employed in a series of relaxed Jacobi updates (as given in Equation \eqref{eqn:relaxed-jacobi}) to accelerate the convergence of the standard Jacobi iteration. One difficulty that arises with the use of these SRJ schemes are issues of floating point error associated with the large range of relaxation factors. For example, applying all of the large relaxation parameters in a given SRJ scheme can potentially lead to overflow. To ameliorate these problems, we follow Yang's suggestion in \cite{Yang14} of cycling through relaxation parameters and successively applying the parameters which will maximize the error reduction at each step for robust convergence in the presence of roundoff.

\section{Developing a data-driven heuristic for selecting SRJ schemes}

The Scheduled Relaxation Jacobi method provides an approach to accelerate Jacobi iteration by applying sets of relaxation factors that can improve convergence. A family of SRJ schemes can be derived, with each scheme corresponding to a different number of distinct relaxation factors $M$ desired in a cycle. Given a linear system of interest, one may seek an SRJ scheme which provides fast convergence. One approach for choosing an optimal scheme for a given linear system is to compute the asymptotic convergence rate associated with many schemes according to Equation \eqref{eqn:convergence-rate} and select the scheme which yields the largest asymptotic convergence rate per iteration. However, for many large scale problems, this computation may be prohibitively expensive as the eigenvalues of the Jacobi iteration matrix may not be easily computable. In this case, the only alternative is to perform extensive experimentation and directly test many schemes to determine an appropriate SRJ scheme which provides fast convergence. We desire an approach to determine a reasonable scheme for convergence without testing multiple schemes or having to compare their convergence rates.


To illustrate this problem further, we show the performance of several SRJ schemes when solving the one-dimensional Poisson equation. The 1D Poisson equation is a prototypical test PDE and is given in Equation \eqref{eqn:poisson}.
\begin{equation}
    -\frac{d^2 u}{dx^2} = f(x) \ , \ x \in [0,1]
    \label{eqn:poisson}
\end{equation}
A finite difference discretization of Equation \eqref{eqn:poisson} on a one-dimensional uniform mesh with Dirichlet boundary conditions leads to a linear system $Ax = b$ where $A$ is a symmetric and tridiagonal matrix with the following entries
\[
A = \frac{1}{\Delta x^2}
\begin{pmatrix}
2 & -1 &  \\
-1 & 2 & -1  \\
 & \ddots & \ddots & \ddots  \\
 &  & -1 & 2 & -1  \\
 &  &  & -1 & 2
\end{pmatrix}
\] 
where $N$ represents the number of degrees of freedom and $\Delta x = \frac{1}{N+1}$. We consider a linear system of size $N = 100$, and show the convergence behavior associated with applying several different SRJ schemes corresponding to different $M$. We set the right hand side vector $b = 1$ for all entries.

\begin{figure}[htbp!]
    \centering
    \includegraphics[width=0.5\textwidth]{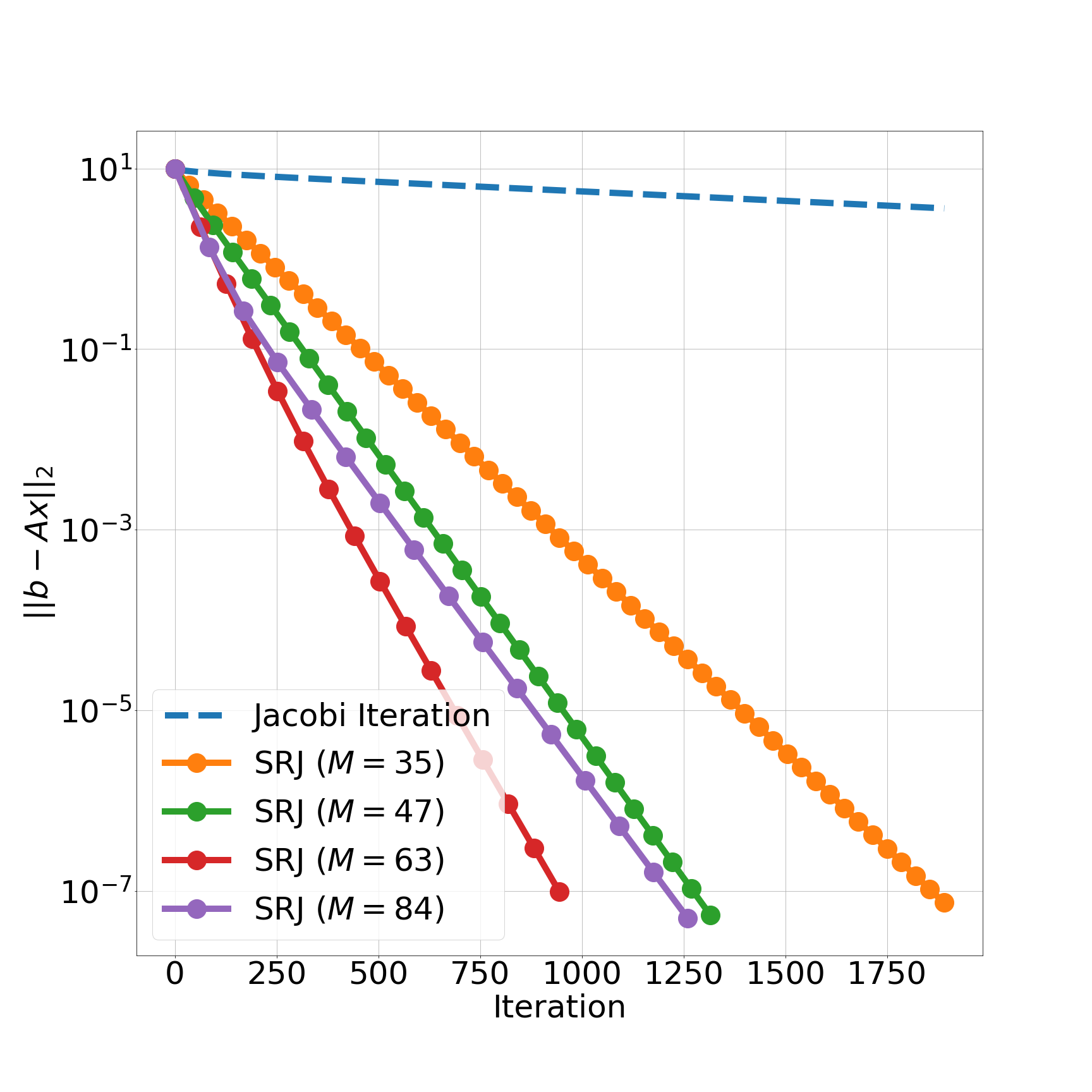}
    \caption{Convergence of several SRJ schemes ($M = 35, 47, 63, 84$) for solving the 1D Poisson problem with $N = 100$. The $M = 63$ scheme exhibits the fastest convergence.}
    \label{fig:srjConvergence_1DPoisson}
\end{figure}

Figure \ref{fig:srjConvergence_1DPoisson} shows the convergence history of the $L_{2}$ residual norm at each SRJ cycle when applying the SRJ schemes associated with several different $M$ in order to solve the linear system arising from discretization of the 1D Poisson equation. All schemes shown illustrate an improvement in convergence relative to the Jacobi iteration. According to Figure \ref{fig:srjConvergence_1DPoisson}, the $M = 63$ scheme provides the most rapid convergence between the four schemes. Our goal is to determine a rapidly converging scheme without needing to test multiple schemes as done here. While one could perform convergence analysis of the four schemes to select the one which provides the fastest asymptotic convergence for this problem according to Equation \eqref{eqn:convergence-rate}, we wish to circumvent this calculation as it may be intractable for much larger problems. 


We propose a data driven approach to determine which SRJ scheme to apply in a given cycle. The approach serves as a surrogate which avoids both computation of asymptotic convergence rates/spectral radii associated with different SRJ schemes to determine the optimal one, as well as brute force experimentation when such analysis is not possible. The key idea is to develop a heuristic that is based on convergence data collected from randomly applying SRJ schemes. At a given SRJ cycle, the heuristic can be used to decide the SRJ scheme at the next cycle which will provide the best overall convergence rate, based on current solution state parameters. To simplify the decision making process, we filter out select schemes corresponding to particular $M$ and designate these as particular scheme levels. The scheme levels are chosen based on their capability to solve stiff systems, which can be characterized by the slope of the amplification polynomial $G_{M}(\lambda)$ at 1. We filter schemes such that there is at least a $50\%$ increase in this value between schemes at adjacent scheme levels. Practically, this filtering process removes schemes which would behave too similarly to one another for stiff problems. The values of $M$ corresponding to a specific scheme level are shown in Appendix B. When selecting an SRJ scheme to use for the next cycle, the solver is restricted to three choices. Given the current scheme level, the solver can either choose to increase the level, keep the same level, or decrease the level and use the corresponding scheme at the next cycle. Our goal is to develop a heuristic that allows us to select the scheme that will provide the best convergence rate at the next cycle.

\begin{figure}[htbp!]
    \centering
    \includegraphics[width=\textwidth]{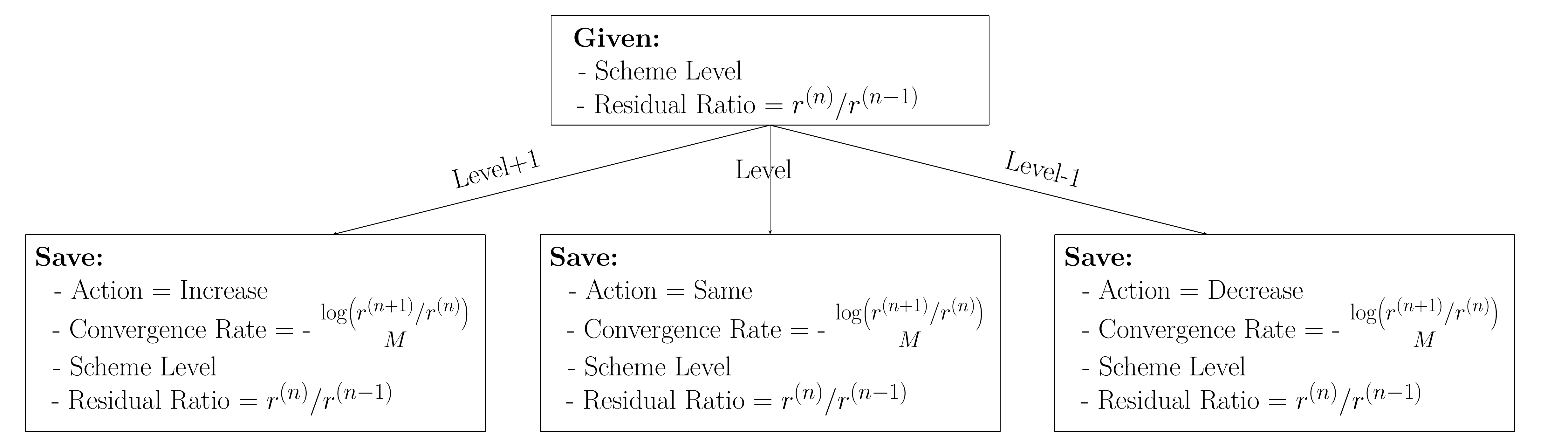}
    \caption{Depiction of data collection procedure. Three data points are collected at each step corresponding to each of the three possible schemes we can select for the next cycle (resulting from increasing, decreasing or keeping the scheme level). Each data point contains four parameters - the action taken, the average convergence rate of SRJ using the new level, the current scheme level, and the residual ratio obtained at the previous cycle. Note that the average convergence rate associated with each action/scheme taken will differ.} 
    \label{fig:srj-data-collection-steps}
\end{figure}

To develop a rule for selecting SRJ schemes, we collect convergence data from the tridiagonal linear system arising from discretization of the 1D Poisson equation as given in Equation \eqref{eqn:poisson} for varying discretizations. We consider matrices of sizes $N = 2,5,10,20,30,40,50,60,70,80,90,100,200,300,400$. We perform a data collection process for each matrix. In each case, the initial solution is a vector of zeros while the right hand side is a vector of ones. We begin by applying the simplest SRJ scheme corresponding to $M = 1$ (scheme level 0) at the initial step. Applying a cycle of SRJ involves a series of $M$ relaxed Jacobi updates given by Equation \eqref{eqn:relaxed-jacobi}. Afterwards, the ratio of the $L_{2}$ residual norm after and before applying the scheme is recorded along with the scheme level that was used. All subsequent data collection steps involve the following procedure (depicted visually in Figure \ref{fig:srj-data-collection-steps}). The three choices for the SRJ scheme to use at the next cycle (corresponding to increasing the level, decreasing the level, or keeping the same level) are employed and the average rate of convergence per iteration associated with using each scheme is computed. We compute the average convergence rate of the iterations comprising an SRJ cycle according to Equation \eqref{eqn:average-convergence-rate}, where $r^{(n)}$ is the $L_{2}$ residual norm at the beginning of the cycle, $r^{(n+1)}$ is the  $L_{2}$ residual norm after performing the SRJ iterations in the cycle, and $M$ is the number of iterations performed in the cycle, equivalent to the number of distinct relaxation parameters in the applied SRJ scheme. Note that our definition of the average convergence rate differs from a geometric mean measure as used in \cite{Yang14} and \cite{Adsuara1} to characterize the per iteration convergence rate over an SRJ cycle.
\begin{equation}
    \text{Average Convergence Rate} = - \frac{\log \left(r^{(n+1)}/r^{(n)} \right)}{M}
    \label{eqn:average-convergence-rate}
\end{equation}
The action taken to modify the scheme level (increase, decrease, same) and the associated convergence rate are recorded. The residual ratio observed at the previous cycle (i.e. ratio of the residual before and after performing SRJ iterations at the previous cycle) is recorded as well. These are the same for the three actions and are considered our state variables. Given the previous residual ratio and current scheme level, we now have the average convergence rate that will be achieved by taking each of the three possible actions at the next cycle. The best scheme to use at the next cycle given our state is now easy to determine - it is the one which provides the largest convergence rate. We save the data associated with each of the three actions as a tuple of data as depicted in Figure \ref{fig:srj-data-collection-steps}. To finish this data collection step, we randomly select one of the three schemes to use and record the residual ratio obtained from utilizing this scheme as well as the level associated with the scheme. This is the given data we store in preparation for the next data collection step, which follows the same procedure as before (determine the convergence rates associated with the three possible actions, record data, and randomly select an action and record the residual ratio and level). It is worth noting that at the first step, we will not be able to decrease the level (due to being at the lowest scheme level) so only two actions (increasing the level or keeping the scheme level) are available.

This illustrates one step of our data collection process, which is shown in Figure \ref{fig:srj-data-collection-steps}. In summary, we collect data regarding our current state (the residual ratio at the previous step and the current scheme level) and the convergence rate associated with each of the three possible schemes (each step provides three data points). Data is collected by performing many trials in which SRJ is used to reduce the residual of the solution below an $L_2$ norm of $10^{-8}$. Each trial consists of many cycles of the SRJ algorithm applied using the data collection procedure described. During our experimentation, we collected roughly one million data points for each size $N$. Given a current solver state which we quantify by the residual ratio observed at the previous cycle and the current level, we seek to determine which action will result in the best convergence rate based on the collected data.

\begin{figure}[htbp!]
    \centering
    \includegraphics[width = 0.5\textwidth]{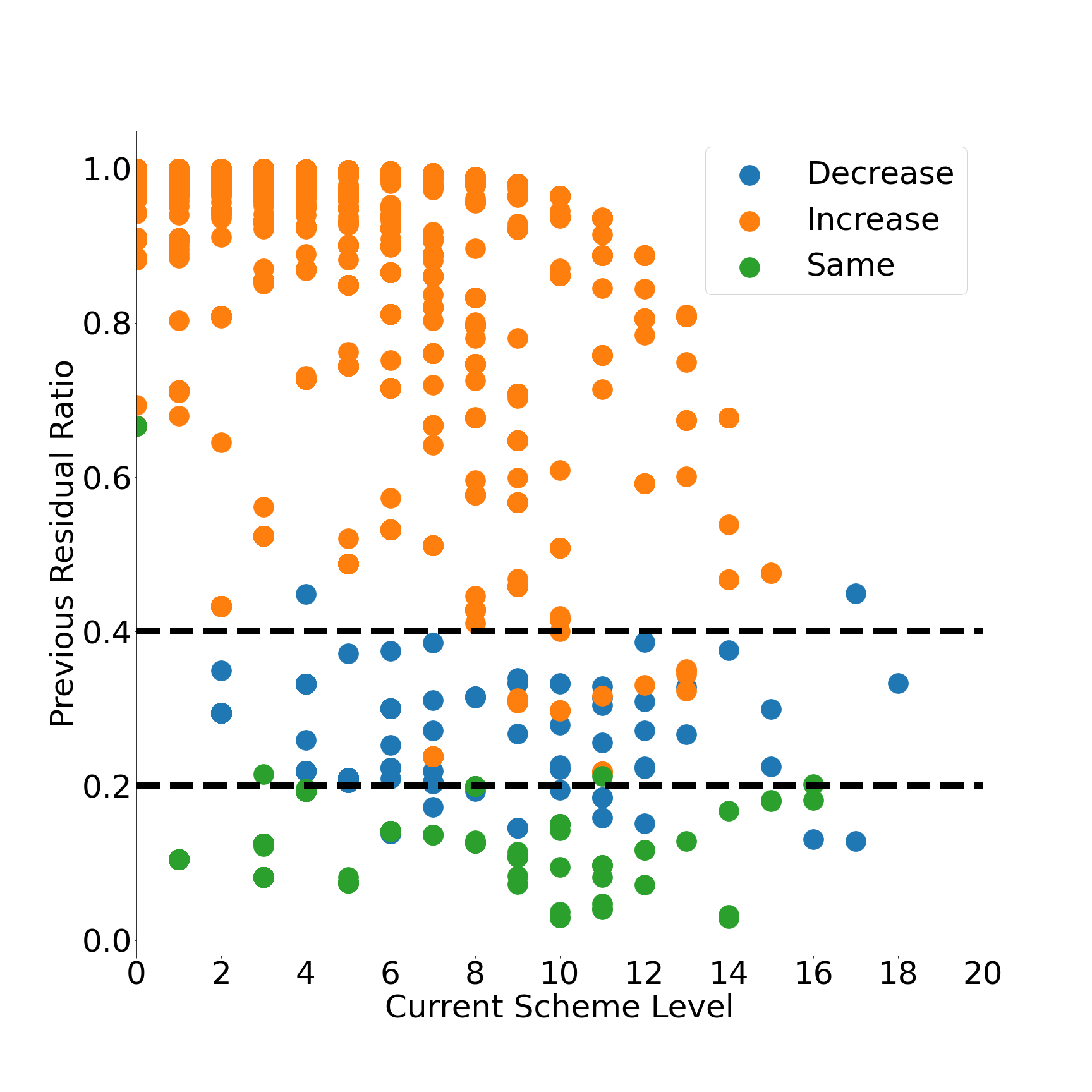}
    \caption{Visualization of the best action to take within the level vs residual ratio space. We develop simple heuristics to determine which action to take based on the previous residual ratio. The current level does not appear to affect the best action.}
    \label{fig:data-clusters}
\end{figure}

The collected data is postprocessed to develop a heuristic for selecting SRJ schemes during the linear solve process. All of the convergence data obtained is organized according to the current scheme level state variable. All data points in each level are sorted in order of increasing previous residual ratio value and grouped into clusters of $N_{\text{set}} = 10000$ data points. Within each cluster, the data points corresponding to increasing the level, decreasing the level and keeping the level are identified, and the mean of the average convergence rates associated with each of the actions is computed along with a $95\%$ confidence interval to ensure we are confident about which action is the best within the cluster. The action corresponding to the highest average convergence rate in each cluster is recorded after ensuring that the confidence interval associated with it does not overlap with the confidence intervals of the other two actions. The mean residual ratio of all data points in the cluster is also recorded. The best action for every cluster of points can be plotted as a single data point in the level and residual ratio space as shown in Figure \ref{fig:data-clusters}.

The data aggregation step described above allows us to visualize how the state variables are related to the most desirable action for fastest convergence. In fact, the region of the level and residual ratio state space where a specific cluster lies appears to have a big influence on the best action to take at the next cycle. In general, when the residual ratio in the previous step is high, the best action for obtaining a good convergence rate at the next cycle is to increase the scheme level. If the residual ratio at the previous step is very low, it is best to continue using the same set of relaxation factors as before. In between these regions, it is best to decrease the level. We develop a simple rule to decide how to select between SRJ schemes. If the previous residual ratio is above 0.4, we should increase the scheme level. If the residual ratio is between 0.2 and 0.4, it is best to decrease the level. When the residual ratio is lower than 0.2, it is best to keep using the same SRJ scheme for good performance. This rule is summarized in Algorithm \ref{alg:heuristic}.

\begin{algorithm}[htbp!]
\caption{Select SRJ Schemes using data driven heuristic}
\label{alg:heuristic}
\begin{algorithmic}
\STATE Given Previous Residual Ratio and Current Scheme Level
\IF{Residual Ratio $>$ 0.4}
\STATE Level = Level + 1
\ELSIF{{Residual Ratio $<$ 0.4} and {Residual Ratio $>$ 0.2}}
\STATE Level = Level - 1
\ELSE
\STATE Level = Level 
\ENDIF
\end{algorithmic}
\end{algorithm}

The data suggests a general and simple heuristic that can be used to determine the best action to take at a given step for selecting SRJ schemes and obtaining robust convergence. The current scheme level being used does not appear to affect the best action to take at the next step, so our rule is independent of this parameter. Our heuristic, summarized in Algorithm \ref{alg:heuristic}, can be used to automatically select SRJ schemes for each problem in order to solve linear systems efficiently. Given the SRJ schemes (which can be derived using Equation \eqref{eqn:srj-relaxation-parameters}) associated with each scheme level shown in Appendix B, the heuristic can be used to select the scheme to apply within the solver at each cycle, without any user intervention.

\section{Numerical Results}

We investigate the convergence behavior of the SRJ method with our scheme selection heuristic developed in Section 3 for solving a variety of linear systems. We begin by exploring the behavior of our solver on linear systems involving the training matrices from which we initially collected convergence data from (specifically, those arising from discretization of the 1D Poisson equation). Afterwards, we explore the convergence behavior on test matrices which are progressively more different from these in-sample training matrices. In particular, these test matrices correspond to the following problems:
\begin{enumerate}
    \item 1D Poisson equation on uniform grid discretized with finite difference method (of different sizes from those in the training set) and random tridiagonal matrices which are symmetric and diagonally dominant
    \item 2D Laplace equation on uniform grid discretized with finite difference method
    \item 3D Poisson equation on uniform grid discretized with finite difference method 
    \item 2D Poisson equation on unstructured grids discretized with finite element method
\end{enumerate}

To assess the efficacy of the scheme selection heuristic with the SRJ method, we compare the convergence behavior of this approach to several other methods. In particular, we compare convergence behavior to that of the standard Jacobi iteration method, as well as the SRJ method with a level selection rule which begins at the lowest scheme level and always increases the level at each subsequent cycle. This approach can be interpreted as a simple brute force heuristic one could use to select schemes without the need to experiment with each scheme, and provides a benchmark approach that our heuristic should outperform. Lastly, we also compare these methods to the most recent SRJ schemes developed by Adsuara et al. in \cite{Adsuara2} which we refer to as CJM (Chebyshev-Jacobi method) schemes to distinguish them from the SRJ schemes developed in this work. The CJM schemes are guaranteed to provide optimal convergence for a given discretization (specified by $N$) and number of relaxation parameters chosen (specified by $M$). We do not expect to outperform this method when an optimal $M$ is chosen for CJM, but rather show that our data driven heuristic provides desirable convergence properties close to this optimal. The main advantage of our approach is the ease of use of the heuristic for automatically selecting appropriate schemes, which provides flexibility for solving different problems. In several cases, we are also able to achieve improved convergence relative to the CJM method, particularly when solving linear systems arising from non-uniform unstructured meshes.

\subsection{In-sample testing}

As a first test of our data based heuristic, we solve the in-sample training linear systems we initially collected data from. The matrices comprising these linear systems are the tridiagonal matrices of size $N = 2,5,10,20,30,40,50,60,$ $70,80,90,100,200,300,400$ corresponding to discretization of the 1D Poisson equation using finite differences on a uniform grid with Dirichlet boundary conditions. For each system, the right hand side is set to a vector of ones of size $N$ and the initial solution vector is the zero vector of the same size. The performance of the SRJ method using our heuristic to choose a scheme for the next cycle is compared to several other methods. 

\begin{figure}[htbp!]
    \centering
    \includegraphics[width=0.5\textwidth]{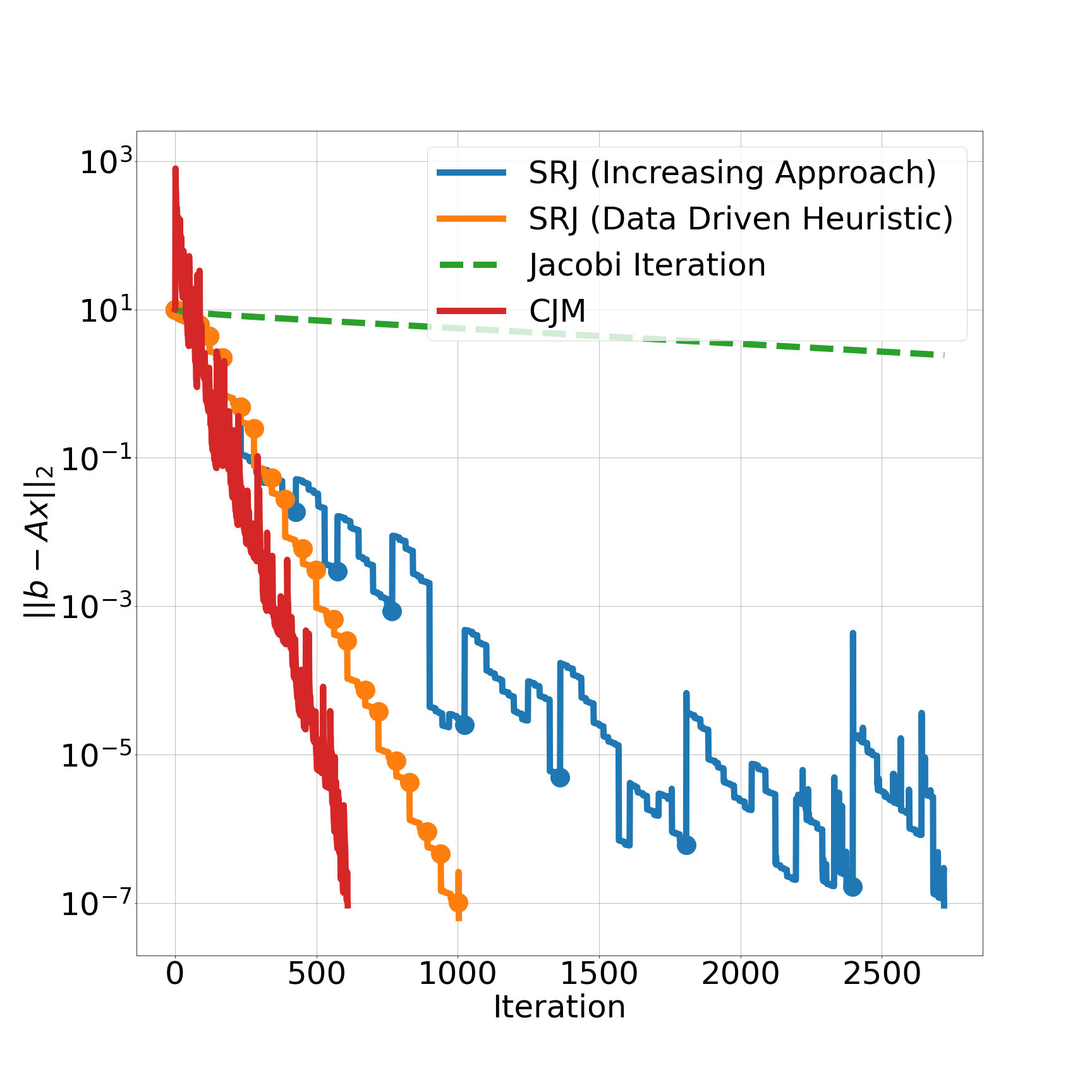}
    \caption{Convergence of SRJ, standard Jacobi, and CJM approaches for the 1D Poisson matrix of size $N = 100$. The CJM solver with optimal $M$ provides the best convergence. However, the SRJ solver with our heuristic works well without requiring experimentation with different schemes, and also outperforms SRJ with a brute force increasing level approach or a standard unrelaxed Jacobi iteration.}
    \label{fig:training-convergence-N100}
\end{figure}

Figure \ref{fig:training-convergence-N100} illustrates a convergence plot for the case of $N = 100$. We show the convergence of the SRJ method with our heuristic as well as an increasing rule for scheme selection, along with the Jacobi iteration, and the CJM scheme of size $M$ which gives the optimal convergence (found by experimentation with CJM schemes of different sizes). Convergence is assumed when the $L_{2}$ norm of the residual $||b-Ax||$ falls below a tolerance value of $10^{-7}$. For this problem, the optimal CJM scheme converges the fastest of the four approaches. SRJ with the heuristic outperforms SRJ with the brute force increasing rule as well as Jacobi iteration. The heuristic approach begins at a scheme level of 0, and increases the scheme level until reaching level 11, and which point it alternates between using the SRJ schemes corresponding to levels 10 ($M = 47$) and 11 ($M = 63$). This approach requires approximately 1000 iterations for convergence. Meanwhile, SRJ with a purely increasing rule requires over 3000 iterations for convergence. The standard Jacobi method without relaxation would require many more iterations. Although SRJ with the heuristic does not outperform the optimal CJM method, it converges reasonably rapidly without requiring us to perform any tuning or experimentation with different SRJ schemes. Additionally, this approach greatly outperforms the brute force approach to scheme selection which is uninformed by any data collection.

\begin{figure}[htbp!]
    \centering
    \includegraphics[width=0.5\textwidth]{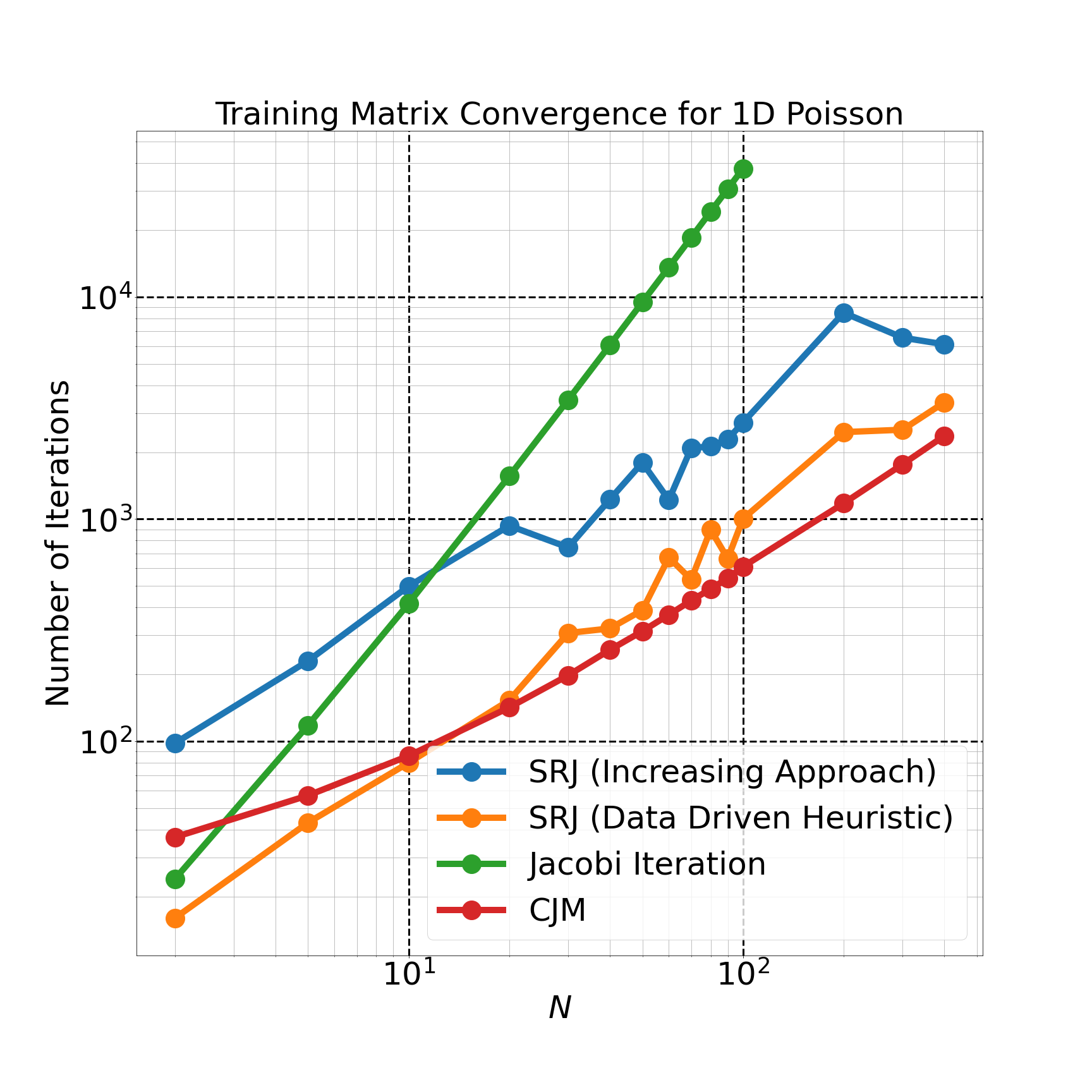}
    \caption{Convergence of 1D Poisson Training Matrices with four different approaches. The CJM method (with optimal $M$) achieves convergence in the fewest number of iterations for all in-sample matrices where $N > 10$. The SRJ method with heuristic converges in a similar number of iterations but does not require testing different schemes for each $N$.}
    \label{fig:1dpoisson-training-convergence}
\end{figure}

This comparison between the four solver approaches is performed for all of the training matrix sizes. In each case, the number of iterations required for each method to converge is recorded. The results are shown in Figure \ref{fig:1dpoisson-training-convergence}. When investigating the CJM method for each $N$, we experiment with many schemes corresponding to different $M$ to determine which provides the fastest convergence. As a result, the CJM method illustrates convergence in the fewest iterations for nearly all $N$ (except when $N < 10$) . However, SRJ with the heuristic for scheme selection shows convergence in a similar number of iterations. In the worst case, this approach requires twice as many iterations to converge, but for many $N$ requires convergence in a similar number of iterations to CJM. The heuristic also outperforms a brute force SRJ approach with a purely level increasing rule, generally requiring at most half the number of iterations as the increasing approach. The Jacobi method works well for small $N$ but the number of iterations increases rapidly as $N$ grows larger. For $N > 100$, Jacobi iteration does not converge in a reasonable number of iterations so we do not record this data in Figure \ref{fig:1dpoisson-training-convergence}. Overall, SRJ with the heuristic performs well on the training matrices and exhibits convergence behavior close to the optimal CJM scheme for each $N$ without requiring any experimentation.

\subsection{Out of sample testing}

Generalizeability is an important aspect of any tool based on data. While the heuristic works well for in-sample matrices, it is important that it generalizes well to arbitrary matrices and results in a scheme selection pattern which yields good convergence for problems which are out of sample. In this section, we explore test problems which are progressively more different from our training problems, including linear systems arising from discretization of higher dimensional PDEs as well as discretization on unstructured domains.

\subsubsection{1D Poisson and Tridiagonal Matrices}
\begin{figure}[htbp!]
    \centering
    \includegraphics[width=0.5\textwidth]{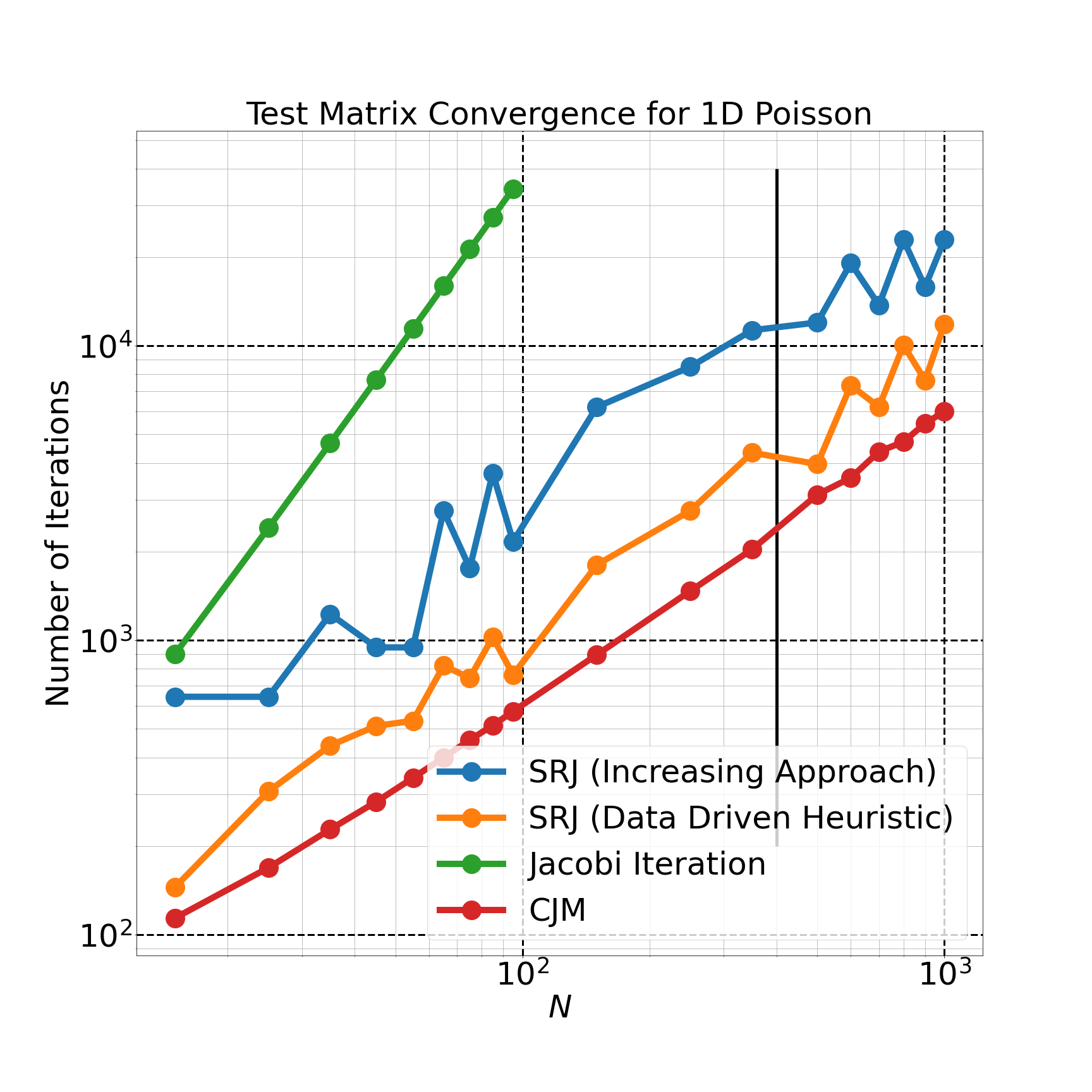}
    \caption{Convergence of out of sample 1D Poisson Matrices. The scaling results are similar to that of the training matrices, suggesting that the heuristic does not overfit to the training matrices but also works well on the test matrices.}
    \label{fig:1dpoisson-test-convergence}
\end{figure}

Linear systems arising from discretization of the 1D Poisson equation (which have sizes differing from our training matrices) are the most similar to those in our training set which are out of sample.  As a result, they provide a good initial test of our heuristic on unseen matrices. We test our approach on 1D Poisson matrices of size $N = 15, 25, 35, 45, 55, 65, 75, 85, 95, 150, 250, 350$ which can be regarded as interpolated samples of our in-sample matrices. Additionally, we test extrapolated samples corresponding to sizes $N = 500, 600, 700, 800, 900, 1000$ (which are larger than the maximum size of our training matrices). As before, the four approaches (SRJ with our data driven heuristic, SRJ with increasing rule, standard Jacobi, and CJM) are used to solve these test linear systems. Figure \ref{fig:1dpoisson-test-convergence} shows the number of iterations required for each method to converge for all of the test $N$. A vertical line is used to demarcate the region where the extrapolated samples have size greater than that of the training samples (the largest training matrix had size $N = 400$). The number of iterations required for convergence on the test matrices scales similarly to the number of iterations required for the training matrices as shown in Figure \ref{fig:1dpoisson-training-convergence}. In particular, the CJM scheme with optimal $M$ selected for each size results in the fastest convergence. However, SRJ with the data driven rule converges in a number of iterations that is close to this optimal. Furthermore, the heuristic works well while the brute force heuristic and standard Jacobi take many more iterations to converge for the matrices considered here. The results indicate that the data based heuristic does not overfit to the training matrices but can still work well on matrices outside of the training set. Additionally, the heuristic does not require any tuning to select an optimal SRJ scheme.

\begin{figure}[htbp!]
    \centering
    \includegraphics[width=0.5\textwidth]{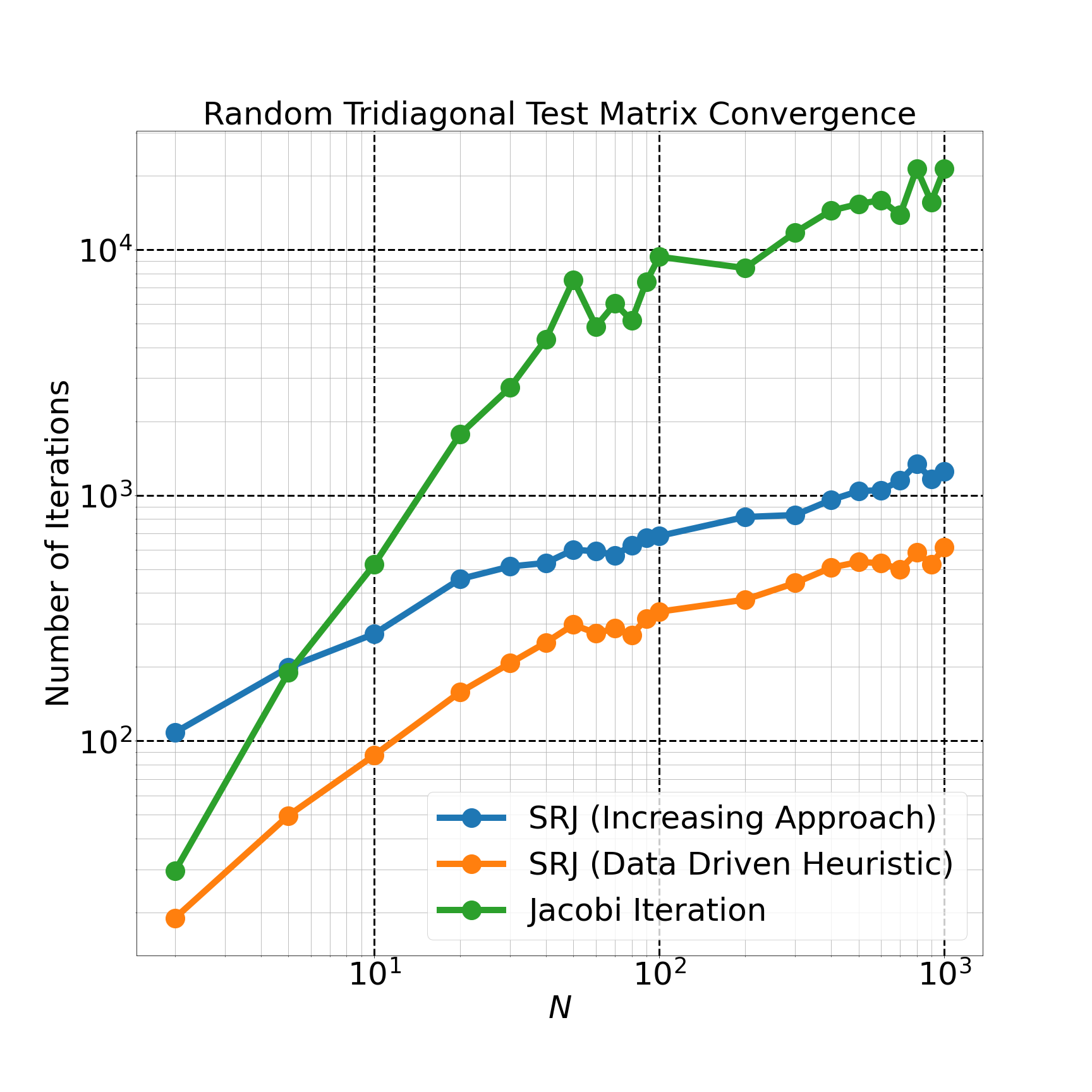}
    \caption{Convergence of random symmetric and diagonally dominant tridiagonal matrices of various sizes. The SRJ method with heuristic provides the best convergence compared to SRJ with the increasing rule and the standard Jacobi method, for all matrix sizes. A comparison with the CJM method is not shown since this method relies on the spatial discretization of a PDE (which is not applicable here).}
    \label{fig:convergence-random-tridiagonal}
\end{figure}

As a more general test of the efficacy of our heuristic, we test our SRJ approach on random symmetric tridiagonal systems which are diagonally dominant, rather than a linear system arising from discretization of a one-dimensional PDE. We constructed random symmetric tridiagonal systems of various sizes, ensuring that they are diagonally dominant so that Jacobi iteration would converge. This was done by creating a random array for the diagonal components as well as a random array for the subdiagonal and superdiagonal (set to the same array to enforce symmetry). For each row that was not diagonally dominant, the diagonal was modified so that its absolute value was equivalent to the sum of the absolute value of the off diagonal entries of that row. The first and last diagonal entries were set to twice the value of the superdiagonal and subdiagonal entry of that row respectively, to resemble the structure of the 1D Poisson equation and to make the linear system easier for Jacobi iteration to solve. Matrices of the following sizes were considered: $N = 2,5,10,20,30,40,50,60,70,80,90,100,200,300,400,500,600,700,800,900,1000$. For each $N$, 20 different tridiagonal matrices which satisfy the symmetric diagonally dominant requirements are constructed and the number of iterations required for standard Jacobi iteration, SRJ with an increasing scheme, and SRJ with our heuristic to converge are recorded. We do not consider the CJM scheme here since the method relies on determining wavenumbers associated with a specific spatial discretization of a PDE (which is not applicable here). The average number of iterations required for each approach to converge is plotted as a function of matrix size $N$ in Figure \ref{fig:convergence-random-tridiagonal}. As before, convergence is assumed when the $L_{2}$ norm of the residual $||b-Ax||$ falls below a tolerance value of $10^{-7}$. For very small $N$, the Jacobi method outperforms the increasing SRJ approach. However, for large $N$, Jacobi iteration generally requires many more iterations (approximately an order of magnitude more) compared to the SRJ approaches. The SRJ approach with the heuristic outperforms the increasing SRJ approach in all cases. For the larger $N$, SRJ with the data based heuristic takes approximately half the number of iterations as SRJ with the increasing level rule. Despite being developed on limited convergence data from the 1D Poisson matrices, the heuristic can be used to obtain relatively good convergence for general symmetric tridiagonal matrices which could be solved by the standard Jacobi iteration. In this case, it is difficult to derive an appropriate CJM scheme since the schemes are designed for linear systems arising from discretization of elliptic PDEs. This presents one advantage of our SRJ schemes. Specifically, they can be applied to solve linear systems which do not necessarily arise from discretization of a PDE, as long as the matrix exhibits symmetry and diagonal dominance.

\subsubsection{2D Laplace equation on uniform domain}

To test a more general class of problems which exhibit behavior further from our in-sample training matrix set, we consider linear systems that arise from discretization of two-dimensional PDEs. We consider the two-dimensional Laplace equation  on the unit square with homogeneous Neumann boundary conditions given by Equations \eqref{eqn:2d-laplace}-\eqref{eqn:2d-laplace-bc-2}. This problem is considered in \cite{Adsuara2} when assessing the performance of the original CJM schemes.
\begin{align}
    \label{eqn:2d-laplace}
    -\nabla^2 u(x,y) \equiv \frac{\partial^2 u}{\partial x^2} + \frac{\partial^2 u}{\partial y^2} &= 0, \ x,y, \in [0,1] \\
    \label{eqn:2d-laplace-bc-1}
    \frac{\partial u}{\partial x} \bigg \lvert_{x = 0} = \frac{\partial u}{\partial x} \bigg \lvert_{x = 1} &= 0, y \in [0,1] \\
    \label{eqn:2d-laplace-bc-2}
    \frac{\partial u}{\partial y} \bigg \lvert_{y = 0} = \frac{\partial u}{\partial y} \bigg \lvert_{y = 1} &= 0, x \in [0,1]
\end{align}

A finite difference discretization is used to discretize equation \eqref{eqn:2d-laplace} on a two-dimensional uniform grid. The Neumann boundary conditions are also implemented using a second order central difference scheme with ghost nodes. The discretization leads to a linear system $Ax = 0$ where the matrix $A$ is sparse and pentadiagonal. For convergence, we consider the maximum relative difference between the current and previous solutions which can be described mathematically as
\begin{equation}
    ||r^{(n)}||_{\infty} = ||x^{(n)} - x^{(n-1)}||_{\infty}
    \label{eqn:relative_residual_solution}
\end{equation}
and is adopted to align with the convergence results presented in \cite{Adsuara2}. Convergence is achieved when the residual falls below a threshold tolerance of $10^{-10}$.
\begin{figure}[htbp!]
    \centering
    \includegraphics[width=0.5\textwidth]{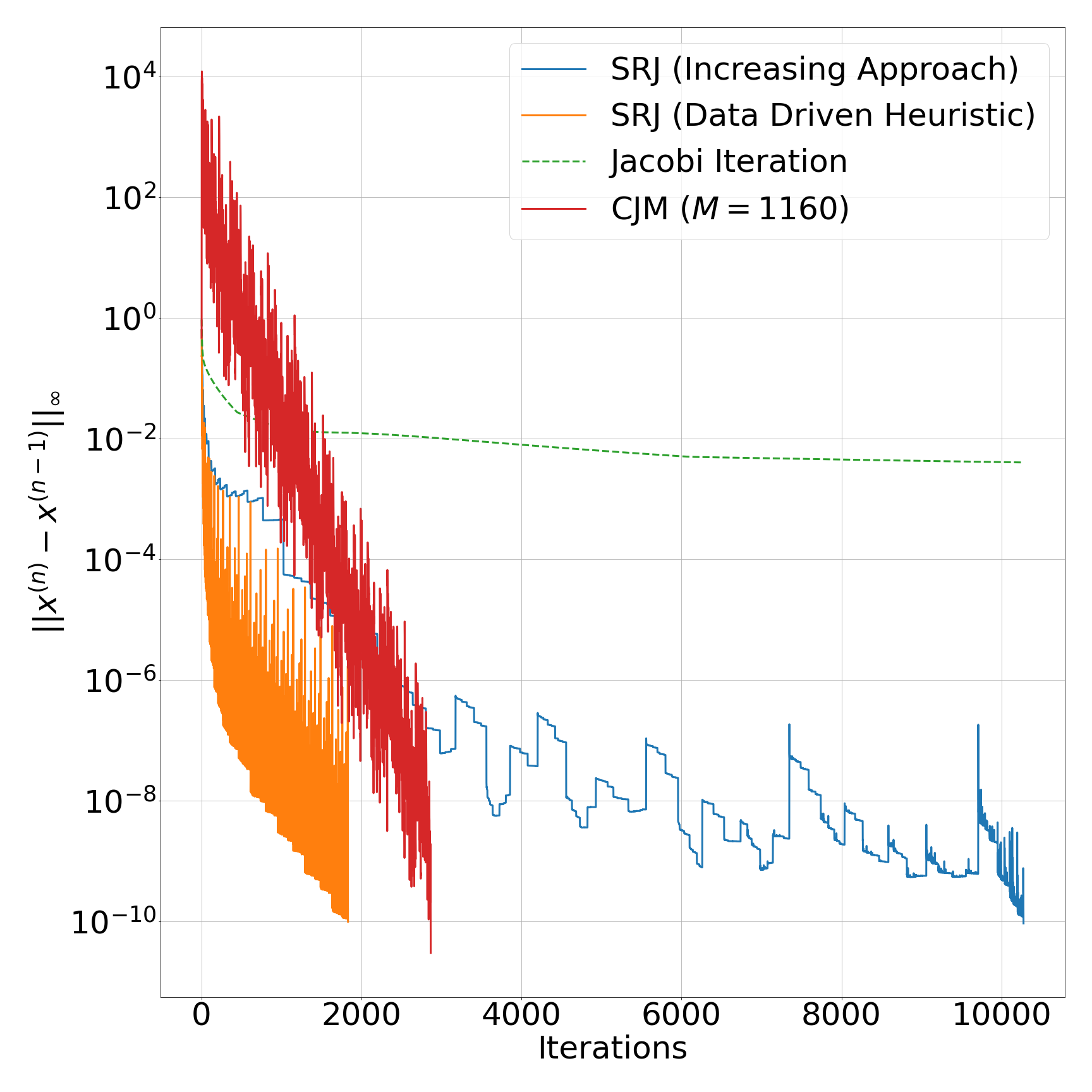}
    \caption{Convergence of the 2D Laplace equation with Neumann boundary conditions for various methods. Our SRJ schemes with the heuristic for scheme selection outperforms all other methods, including the CJM method.}
    \label{fig:2D_laplace_convergence_N256}
\end{figure}
We consider a domain size with 256 $\times$ 256 interior degrees of freedom, so that the linear system has $N = 256^2$ unknowns. The initial solution is set to a random vector, and the residual given by Equation \eqref{eqn:relative_residual_solution} is tracked at every iteration. Figure \ref{fig:2D_laplace_convergence_N256} shows the convergence behavior of our four different methods (SRJ with the data driven heuristic, SRJ with the increasing rule for scheme selection, Jacobi iteration, and the CJM method where $M = 1160$ as used in \cite{Adsuara2}). In this example, the SRJ method with heuristic converges in the fewest iterations, followed by the CJM method. The increasing rule requires many more iterations to converge, and appears to demonstrate worsening convergence behavior as the iterations progress. Finally, the Jacobi iteration stagnates and is unable to converge in a reasonable number of iterations. The results illustrate that there may be an advantage to applying different schemes at different cycles of the algorithm (as is done by SRJ with our heuristic) compared to applying a single scheme for all cycles (as is done by the CJM method). This flexibility allows our SRJ schemes to converge faster than the CJM scheme used here.

To highlight the behavior of the methods for different discretizations, we illustrate the scaling behavior of the number of iterations required for convergence for each method as a function of $N$. In particular, we explore discretizations with the number of DOFs $N = 32^{2}, 64^{2}, 128^{2}, 256^{2}, 512^{2}, 1024^{2}$. We use the CJM scheme corresponding to $M = 3000$ for all discretizations, as this CJM scheme was explored in \cite{Adsuara2} for this benchmark problem. Since the convergence behavior could be sensitive to the initial random vector solution, we perform 10 trials for each discretization and measure the average number of iterations required for convergence for each method. These scaling results are shown below in Figure \ref{fig:2D_laplace_convergence_scaling}. We omit the Jacobi iteration in this scaling study as the method is unable to reach a residual value below the tolerance threshold for this problem. 

\begin{figure}[htbp!]
    \centering
    \includegraphics[width=0.5\textwidth]{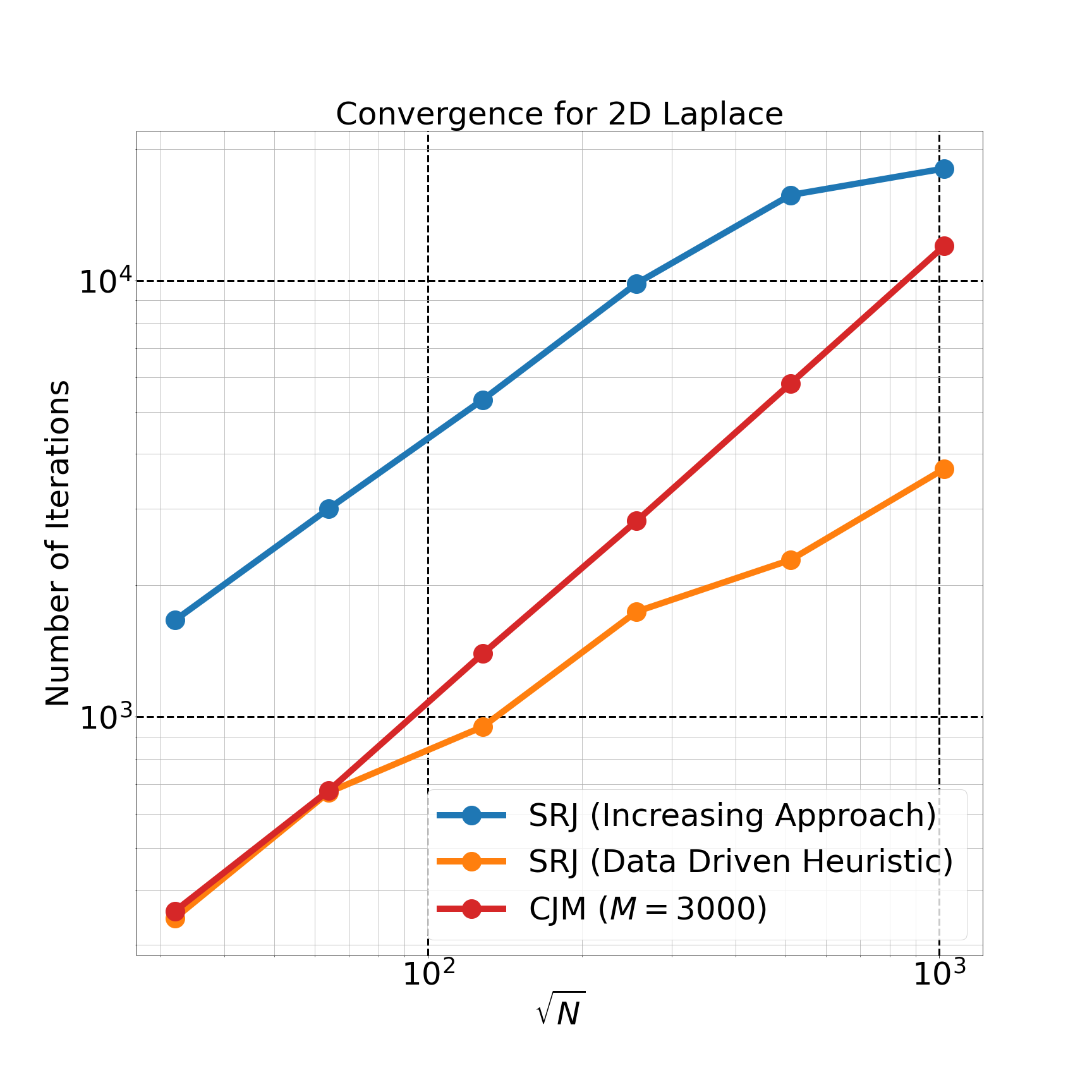}
    \caption{Scaling behavior of the SRJ method with the data driven heuristic and increasing rule, and of the CJM method with $M = 3000$. For smaller discretizations, the SRJ method with heuristic and CJM converge in similar number of iterations. For finer grids, the SRJ method with heuristic converges in fewer iterations.}
    \label{fig:2D_laplace_convergence_scaling}
\end{figure}

For all $N$, the SRJ method with the increasing rule takes the largest number of iterations to converge. For the smaller values of $N = 32^{2}, 64^{2}$, the SRJ method with heuristic and CJM method illustrate convergence in a similar number of iterations. For larger $N$, the number of CJM iterations required for convergence continues to scale with the domain discretization along a single direction (i.e. $\sqrt{N}$). However, the SRJ method with heuristic scales at a slower rate and converges in considerably fewer iterations in this case. It is likely that one could determine an $M$ value for the larger domain sizes such that the CJM scheme converges in fewer iterations. However, determining a good $M$ may be nontrivial and require experimentation, which can be prohibitively expensive for this two-dimensional problem. Our data driven heuristic automatically selects appropriate schemes for convergence without requiring this kind of experimentation by the user. In summary, the SRJ method with heuristic outperforms the CJM method for this 2D Laplace test problem.

\subsubsection{3D Poisson equation on uniform domain}
\begin{figure}[htbp!]
    \centering
    \begin{subfigure}[b]{0.45\textwidth}
        \centering
        \includegraphics[width=\textwidth]{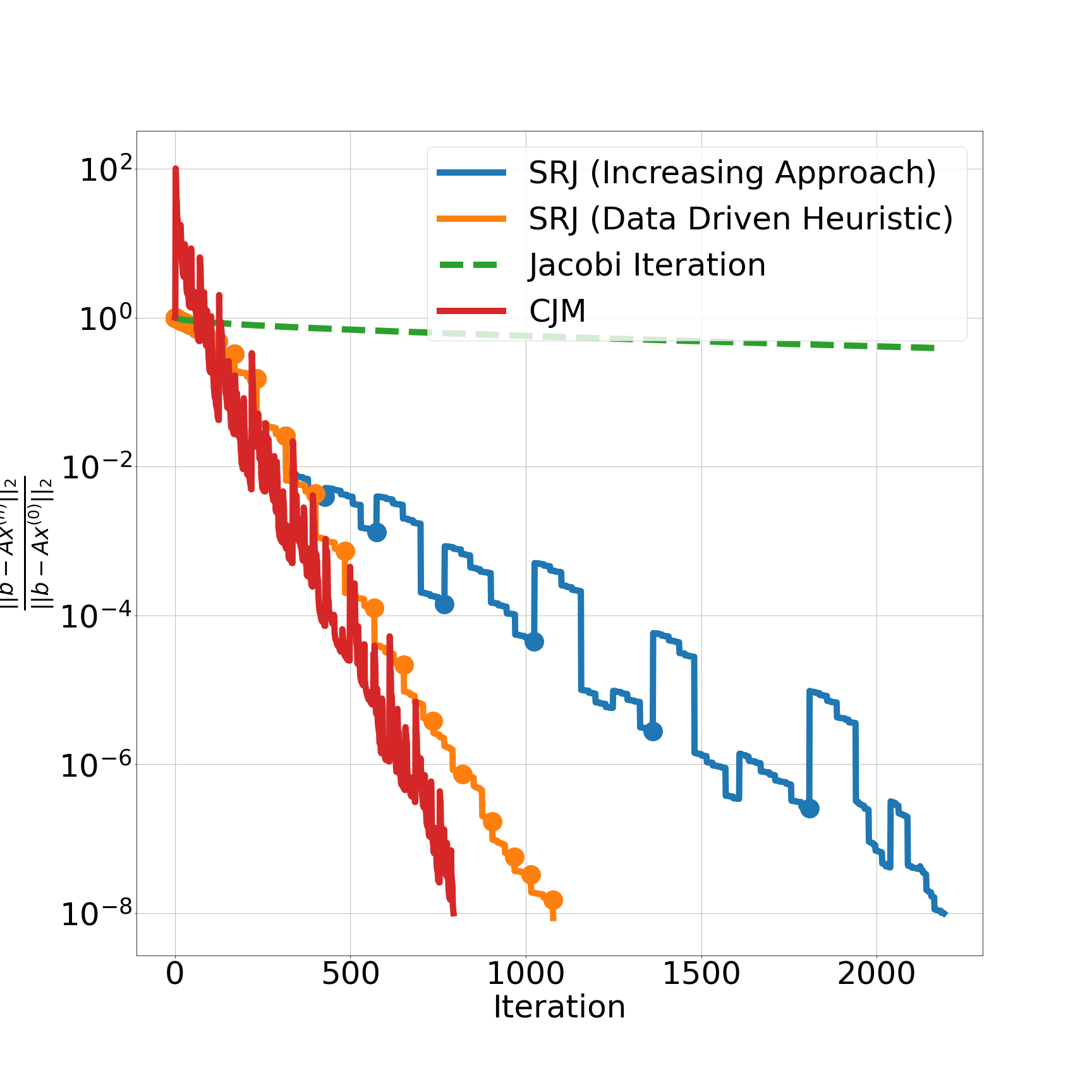}
        \caption{$N = 128^{3}$}
        \label{fig:3D_poisson_convergence_N128}
    \end{subfigure}
    \hfill
    \centering
    \begin{subfigure}[b]{0.45\textwidth}
        \centering
        \includegraphics[width=\textwidth]{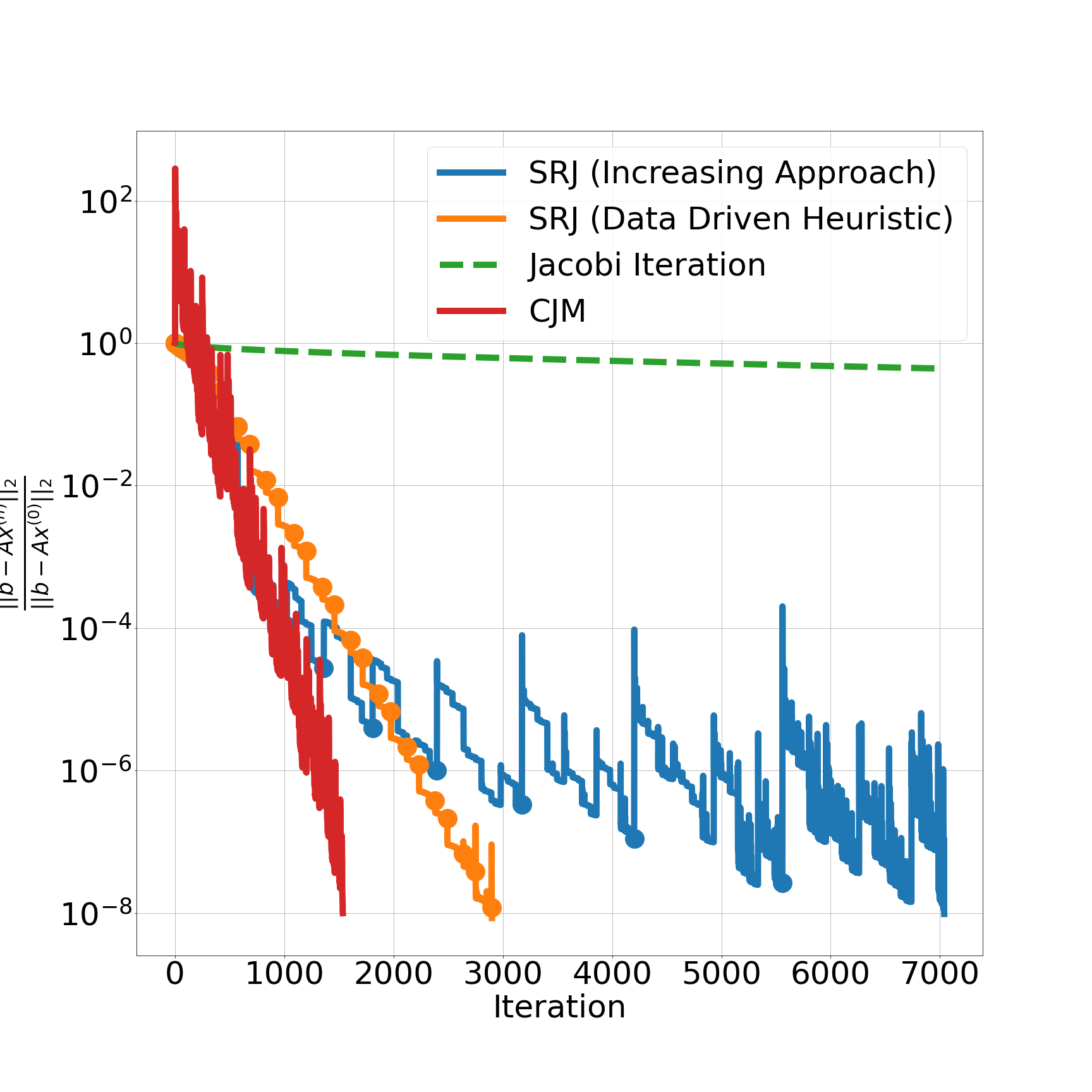}
        \caption{$N = 256^{3}$}
        \label{fig:3D_poisson_convergence_N256}
    \end{subfigure}
    \caption{Convergence behavior of the four approaches in solving the 3D Poisson equation on a uniform cube domain with $N = 128^{3}$ and $N = 256^{3}$ interior DOFs and Dirichlet boundary conditions. CJM with the optimal $M$ converges in the fewest iterations, followed by SRJ with the data driven heuristic, in both cases.}
    \label{fig:3dpoisson-convergence}
\end{figure}

As a further departure from our training linear systems, we consider the 3D Poisson equation on a uniform cubic domain which is given by Equation \eqref{eqn:3d-poisson}. 
\begin{equation}
    -\nabla^2 u(x,y) \equiv -\left( \frac{\partial^2 u}{\partial x^2} + \frac{\partial^2 u}{\partial y^2} + \frac{\partial^2 u}{\partial z^2} \right) = f(x,y,z) \ , \ x,y,z \in [0,1]
    \label{eqn:3d-poisson}
\end{equation}
A finite difference discretization is used to discretize equation \eqref{eqn:3d-poisson} into a set of linear equations, resulting in a linear system $Ax = b$ where the matrix $A$ has up to seven entries per row. Dirichlet boundary conditions are assumed on all edges. We consider two discretizations of Equation \eqref{eqn:3d-poisson}, with 128 $\times$ 128 $\times$ 128 and 256 $\times$ 256 $\times$ 256 interior grid points, resulting in linear systems with $N = 128^{3}$ and $N = 256^{3}$ unknowns respectively. The forcing function is set to $f(x,y,z) = 1$. Figures \ref{fig:3D_poisson_convergence_N128} and \ref{fig:3D_poisson_convergence_N256} illustrate the convergence history for our four different solver methods for solving the linear systems with $N = 128^{3}$ and $N = 256^{3}$ unknowns. The residual corresponding to the end of a cycle is also illustrated with circle markers for both SRJ approaches. Convergence is judged based on the relative $L_{2}$ residual norm which is defined by equation \eqref{eqn:relative-residual} (where $x^{(0)}$ denotes the initial solution and $x^{(n)}$ denotes the current solution). We assume a tolerance threshold value of $10^{-8}$.
\begin{equation}
    \text{Relative Residual} = \frac{||b-Ax^{(n)}||_{2}}{||b-Ax^{(0)}||_{2}}
    \label{eqn:relative-residual}
\end{equation}

In these cases, the CJM method with optimal $M$ (found by experimenting with different $M$ for the analogous 1D Poisson problem and applying the associated scheme to the 3D problem) converges in the fewest number of iterations, followed by the SRJ method with heuristic for scheme selection. For $N = 128^{3}$, the SRJ method with heuristic takes approximately 1000 iterations for convergence while CJM takes approximately 750 iterations. For $N = 256^{3}$, the SRJ method takes approximately 3000 iterations for convergence while CJM takes approximately 1500 iterations. The purely increasing scheme level approach initially demonstrates favorable convergence in both cases, following the convergence behavior of the SRJ with heuristic and CJM methods. However, the convergence of the approach worsens as the iterations progress. The Jacobi iteration would take many more iterations to achieve the desired tolerance. 

\begin{figure}[htbp!]
    \centering
    \includegraphics[width=0.5\textwidth]{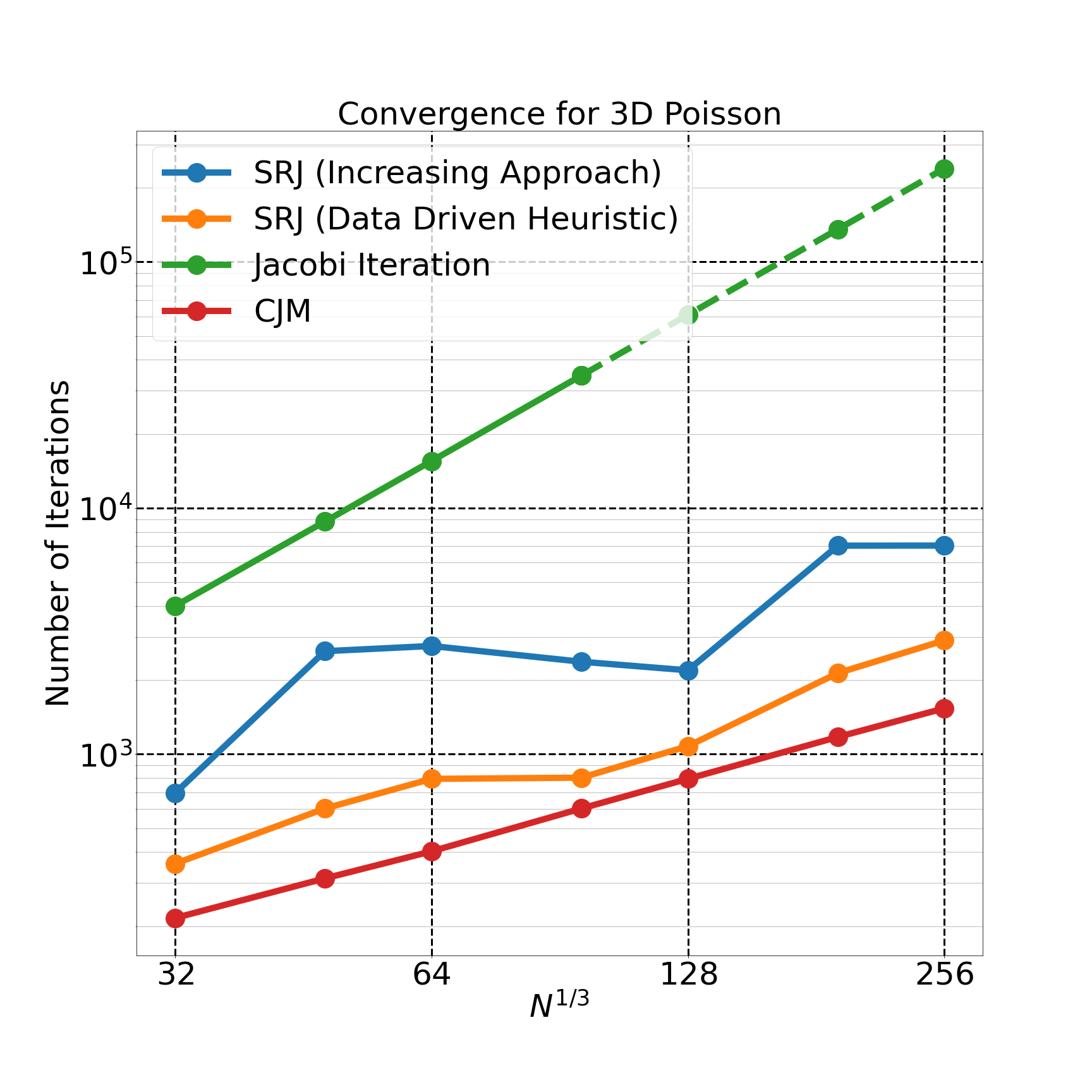}
    \caption{Scaling behavior of the number of iterations required the SRJ method and Jacobi iteration to converge for the 2D Poisson equation discretized on various grid sizes. The SRJ method + heuristic converges much quicker than the other approaches.} 
    \label{fig:scaling-2D-poisson}
\end{figure}

To study the scaling behavior of the four approaches, we consider finite difference discretizations of Equation \eqref{eqn:3d-poisson} involving the following number of DOFs: $N = 32^{3}, 48^{3}, 64^{3}, 96^{3}, 128^{3}, 192^{3}, 256^{3}$. We record the number of iterations required for each approach to converge to a relative $L_{2}$ residual norm below $10^{-8}$. These scaling results are shown in Figure \ref{fig:scaling-2D-poisson} as a function of the number of interior DOFs along a single dimension (i.e. $N^{1/3}$). For Jacobi, we extrapolate the number of Jacobi iterations required for convergence for cases larger than $N = 96^{3}$ (represented as a dashed line). The optimal CJM scheme is chosen by exploring the convergence of different schemes on the analogous 1D Poisson problem. For all problem sizes, the optimal CJM method illustrates the fastest convergence. SRJ with the data driven heuristic takes at most twice as many iterations in the worst case scenario. However, the number of iterations are in some cases only slightly higher, while allowing us to avoid the need to experiment with different schemes for each problem size. The increasing rule illustrates much slower convergence. Additionally, the Jacobi iteration takes many more iterations to converge.

A parameter of interest is the speedup the SRJ method with heuristic offers compared to the standard Jacobi iteration. We define the speedup as the ratio of the number of iterations required for convergence between the two methods. Table \ref{tab:jacobi-vs-srj-3D-speedup} summarizes the speedup obtained by the SRJ method with heuristic as compared to the Jacobi iteration for the different problem sizes. According to Table \ref{tab:jacobi-vs-srj-3D-speedup}, the speedup grows as the problem size grows larger. This could also be observed from Figure \ref{fig:scaling-2D-poisson}, which shows that the number of iterations required for SRJ with the heuristic to converge scales at a slower rate compared to the number of iterations for Jacobi iteration to converge. For the smallest discretization with $32 \times 32 \times 32$ interior DOFs, SRJ with the heuristic provides a roughly 11 times speedup. However, for the finest discretization which employs $256 \times 256 \times 256$ interior DOFs, we obtain a greater than 80 times speedup. We expect these speedups to be even higher for finer discretizations. While the improvement may be slightly smaller than that expected from the optimal CJM scheme, the heuristic approach avoids the need to experiment to find an optimal scheme. In summary, using the heuristic provides a large speedup over the standard Jacobi iteration, and this speedup increases when solving larger linear systems.

\begin{table}[htbp!]
    \centering
    \begin{tabular}{|c|c|c|c|c|c|c|c|c|c|c|c|}
    \hline
    $N$ & $32^{3}$ & $48^{3}$ & $64^{3}$ & $96^{3}$ & $128^{3}$ & $192^{3}$ & $256^{3}$ \\
    \hline
    Approximate Speedup & 11 & 15 & 20 & 43 & 57 & 64 & 83 \\
    \hline
    \end{tabular}
    \caption{Approximate speedup (i.e. ratio of the number of Jacobi iterations to SRJ + Heuristic iterations required for convergence) observed for different discretizations investigated for the 3D Poisson equation. The speedup increases as the number of DOFs increases.}
    \label{tab:jacobi-vs-srj-3D-speedup}
\end{table}

The data driven heuristic, although based solely on convergence data collected from applying SRJ to the 1D Poisson equation, appears to generalize well when solving higher dimensional PDEs. Until now, we have only considered linear systems which arise from discretization of PDEs on uniform domains. In the next section, we show the performance of our linear solver methods on problems arising from finite element discretization of PDEs on several unstructured domains, which may be more representative of the general problems practitioners are interested in solving for large scale high performance simulations.

\subsubsection{2D Poisson equation on unstructured meshes}
\label{sec:srj-convergence-pdes}

In this section, we explore the performance of the SRJ method for solving a PDE on unstructured domains. Specifically, we solve the 2D Poisson equation on several unstructured meshes to determine if the SRJ method with the heuristic for scheme selection can still provide relatively good improvements in convergence when solving more general unstructured PDEs. We also wish to compare the performance of our SRJ schemes to the CJM method, which to our knowledge, has not yet been tested on linear systems arising from discretization of PDEs on unstructured meshes.

\begin{figure}[htbp!]
    \centering
    \begin{subfigure}[b]{0.3\textwidth}
        \includegraphics[width=\textwidth]{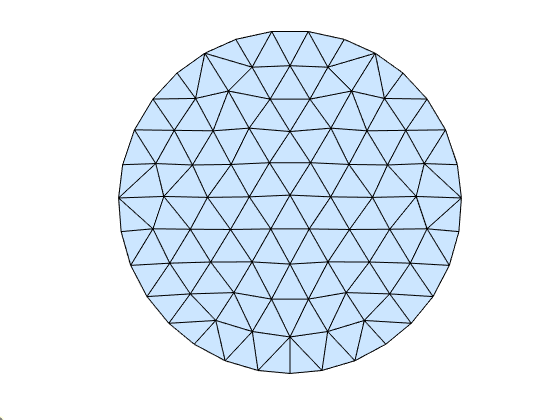}
        \caption{Circle Mesh}
    \end{subfigure}
    \hfill
    \begin{subfigure}[b]{0.3\textwidth}
        \includegraphics[width=\textwidth]{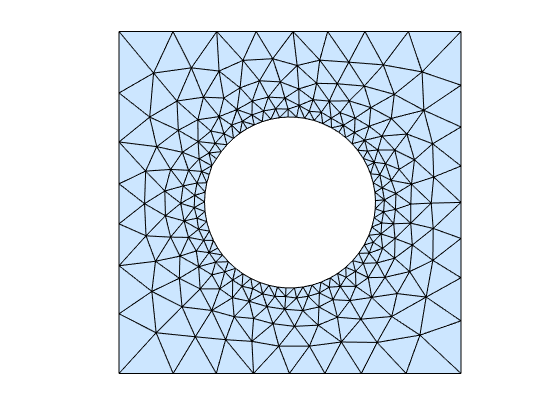}
        \caption{Plate with Hole Mesh}
    \end{subfigure}
    \hfill
    \begin{subfigure}[b]{0.3\textwidth}
        \includegraphics[width=\textwidth]{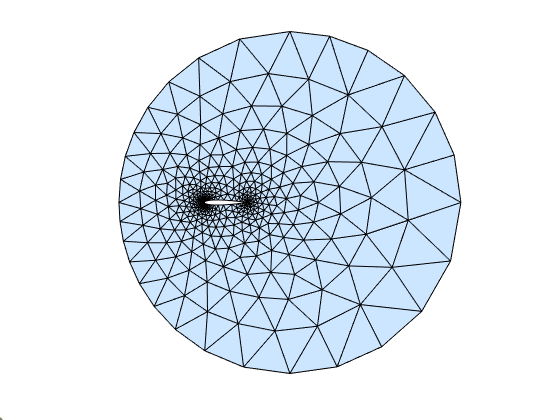}
        \caption{Airfoil Mesh}
    \end{subfigure}
    \\
    \centering
    \begin{subfigure}[b]{0.3\textwidth}
        \vspace{0.1in}
        \includegraphics[width=\textwidth]{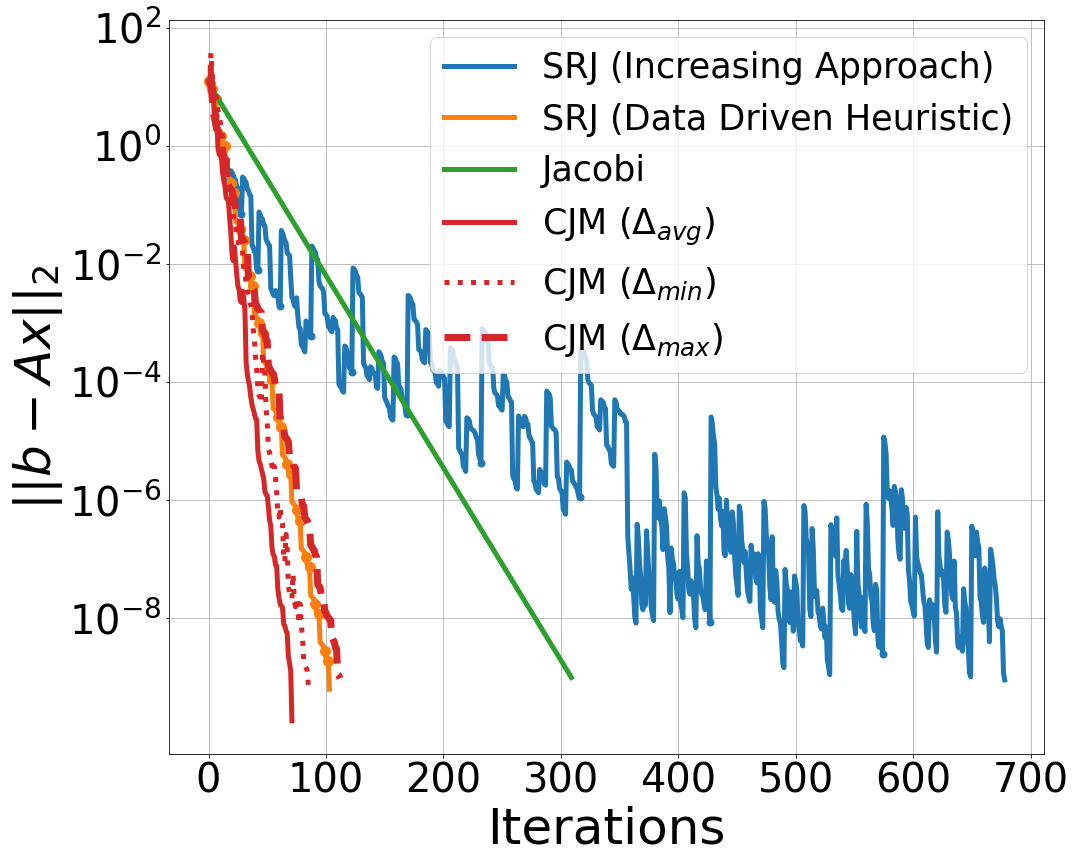}
        \caption{Convergence on Circle Mesh}
    \end{subfigure}
    \hfill
    \begin{subfigure}[b]{0.3\textwidth}
        \vspace{0.1in}
        \includegraphics[width=\textwidth]{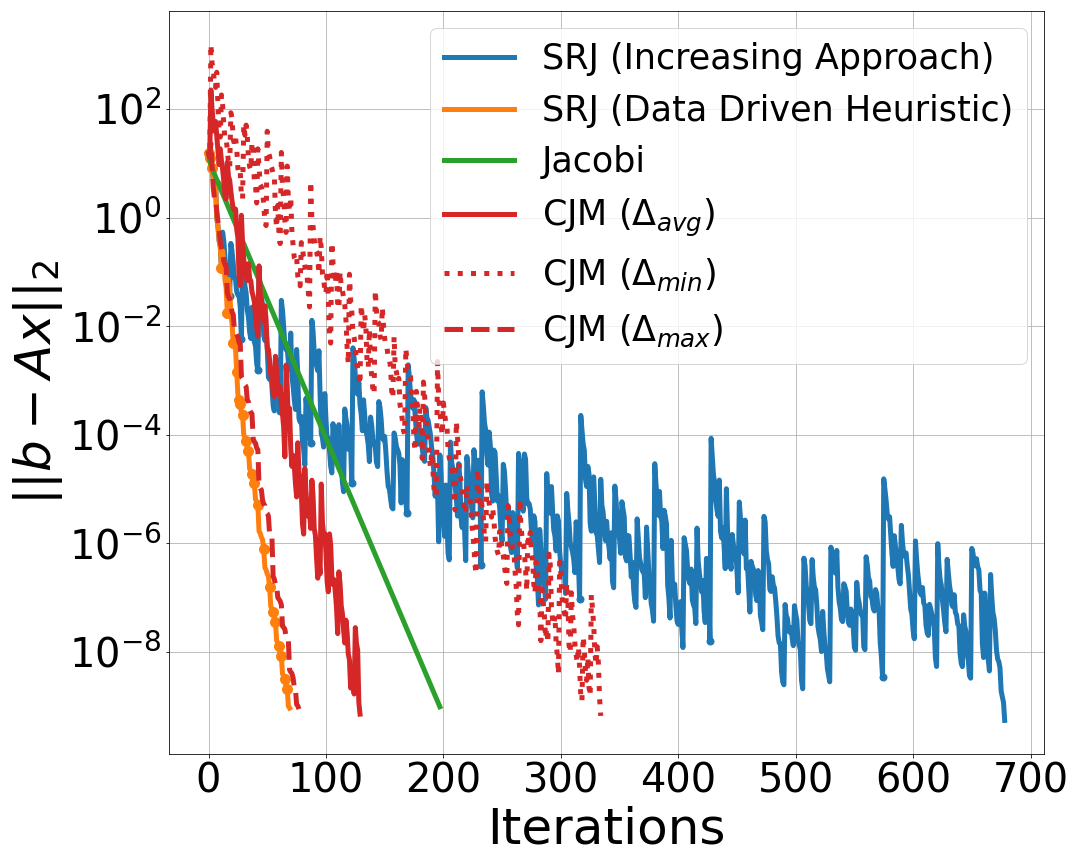}
        \caption{Convergence on Plate with Hole}
    \end{subfigure}
    \hfill
    \begin{subfigure}[b]{0.3\textwidth}
        \vspace{0.1in}
        \includegraphics[width=\textwidth]{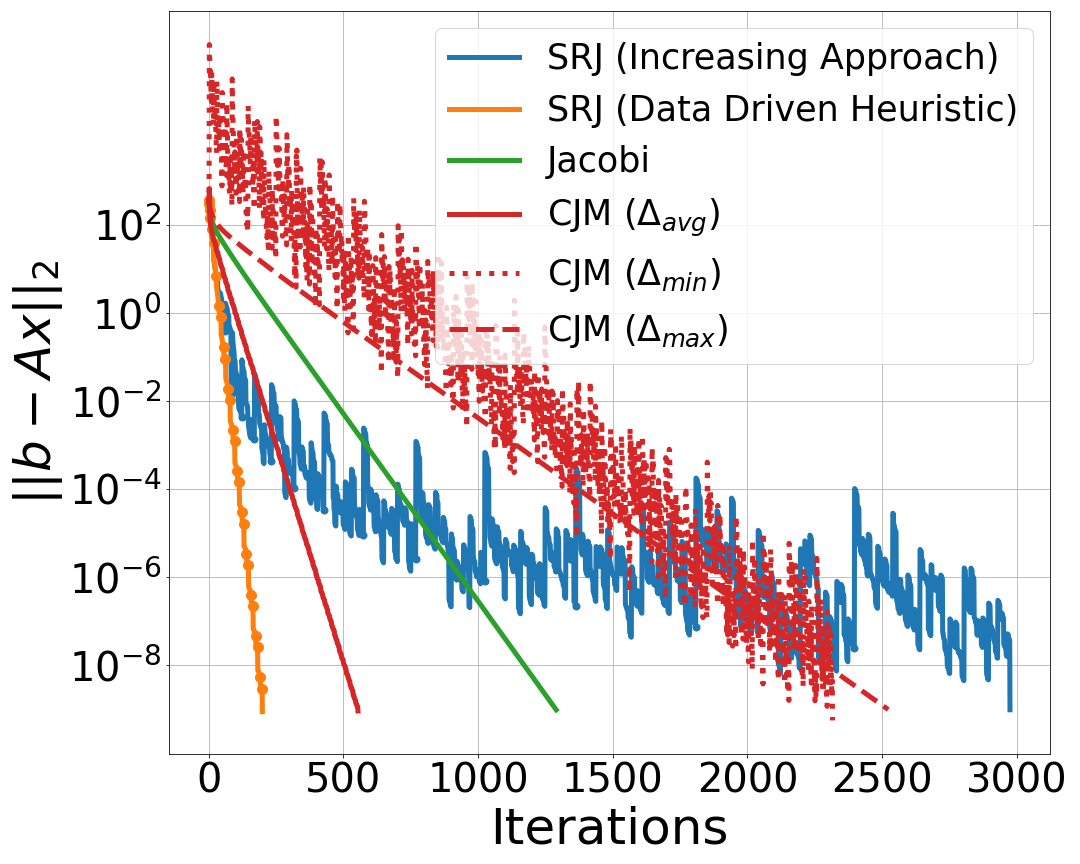}
        \caption{Convergence on Airfoil Mesh}
    \end{subfigure}
    \caption{We solve the 2D Poisson equation on the three different unstructured finite element meshes. The convergence of Jacobi iteration is shown along with the SRJ method (using both the data driven heuristic and a purely increasing rule) and the CJM method (with schemes based on the minimum, maximum, and average grid spacing in the mesh) for each case. While CJM can outperform SRJ on the circle mesh, SRJ with the heuristic gives the best convergence on the plate with hole and airfoil meshes which exhibit a large degree of nonuniformity.}
    \label{fig:sample-meshes}
\end{figure}

We consider three distinct unstructured meshes, corresponding to a circle, plate with hole, and an airfoil geometry. Figure \ref{fig:sample-meshes} illustrates these meshes, which were generated using the MATLAB package \texttt{distmesh} \cite{persson2004simple, persson2005mesh}. We apply the finite element method in order to discretize the 2D Poisson equation into a linear system of equations on each mesh (assuming homogenous Dirichet boundary conditions on boundary nodes) which can be solved to obtain the solution values at the nodes of the mesh. The linear systems obtained using the finite element method on the unstructured meshes are no longer simply tridiagonal or pentadiagonal as before. However, they still exhibit sparsity. The three resulting linear systems are each solved using the Jacobi iterative method, the SRJ method with both data driven heuristic and increasing rule for the scheme level selection, as well as the CJM method. For the CJM method, different schemes can be derived based on a length scale $L$ (here we set $L = 2$ corresponding to the diameter of the circle and airfoil meshes and edge length of the plate with hole mesh) and discretization characterized by some grid spacing $\Delta x$. For a uniform mesh, the grid spacing $\Delta x$ is well defined based on the domain length and number of grid points. However, in the nonuniform case, we instead use three different sets of CJM schemes for each mesh, which are derived based on the minimum, maximum and average grid spacing in the mesh. 

Figure \ref{fig:sample-meshes} shows a visualization of each of the three unstructured meshes we consider in our study, followed by convergence plots of the standard $L_{2}$ residual norm when applying each linear solver method to the linear system associated with each mesh. We assume convergence is achieved at an $L_{2}$ residual norm of $10^{-9}$. For the CJM method, we estimate an $M$ value which is large enough such that convergence could be achieved in a single CJM cycle, based on Equation (14) in \cite{Adsuara2}. The methods illustrate different behavior on different meshes. For the circle and airfoil meshes, the CJM schemes associated with the average length scale gives the best convergence of the three possible CJM schemes, while the CJM scheme associated with the maximum grid spacing gives the best convergence of the three possibilities for the plate with hole mesh. For the circle mesh, the CJM method outperforms SRJ with our heuristic. For the plate with hole and airfoil meshes; however, the SRJ method outperforms all other methods. The speedup achieved using SRJ is particularly prominent for the airfoil mesh, where we achieve convergence in only 199 iterations using SRJ with the heuristic for scheme selection, but require 554 iterations with CJM. This is likely due to the nonuniform nature of the mesh, which makes it difficult to choose an appropriate length scale from which to derive an efficient CJM scheme for this problem. This suggests that the CJM schemes may have difficulty dealing with problems with large degrees of nonuniformity. However, the SRJ schemes appear to work well in this case.

\begin{figure}[htbp!]
    \centering
     \includegraphics[width=0.5\textwidth]{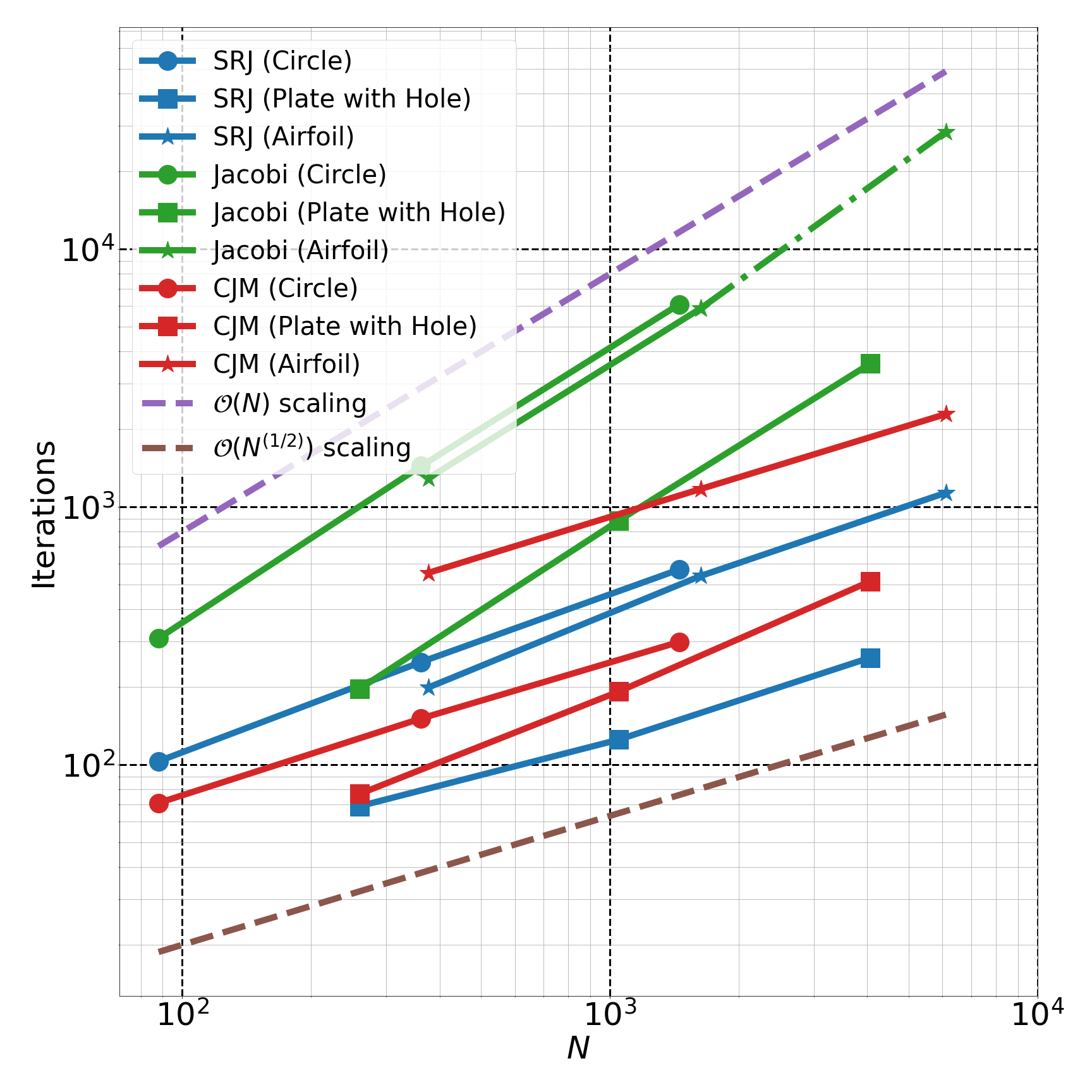}
    \caption{Jacobi vs SRJ scaling with mesh refinement}
    \label{fig:srj_scaling_unstructured}
    \caption{Scaling behavior of the number of iterations required for convergence for Jacobi, SRJ, and CJM for all three geometries. The Jacobi iterative method scales with $\mathcal{O}(N)$ whereas the SRJ method with our heuristic (as well as CJM) scales with roughly $\mathcal{O}(N^{\frac{1}{2}})$. This scaling results in SRJ further outperforming Jacobi iteration as the linear system size is increased, as illustrated by the speedup results in Table \ref{tab:scaling}.}
    \label{fig:srj-jac-scaling}
\end{figure}

We investigate how the number of iterations required for convergence scales as the meshes are refined. For each geometry, we perform two additional levels of refinement (by changing a characteristic length scale parameter within \texttt{distmesh} during mesh generation). This results in three different meshes of various refinement levels for each geometry. The number of iterations required to achieve an $L_{2}$ residual norm below a tolerance level of $10^{-9}$ using the SRJ method with heuristic, Jacobi iteration, and CJM are recorded in each case. For each geometry and linear solver method, the scaling behavior of the number of iterations required for convergence is plotted in Figure \ref{fig:srj_scaling_unstructured} as a function of linear system dimension $N$ (which is larger as the mesh is refined). The CJM results for each geometry correspond to applying the scheme which exhibited the best convergence for each mesh (i.e. for the circle and airfoil meshes, the CJM scheme corresponds to the average grid spacing, while for the plate with hole mesh, the CJM scheme corresponds to the maximum grid spacing). For the finest airfoil mesh, the number of iterations required for Jacobi iteration to converge is extrapolated (as shown with a dash dotted line). In general, we observe that the Jacobi iteration scaling on all three meshes is linear with the number of DOFs $N$. However, the number of iterations required for the SRJ and CJM methods to converge for all meshes scales approximately with the square root of the number of DOFs (i.e. $\mathcal{O}(N^{1/2})$). This implies that the SRJ method with heuristic (and CJM) will further outperform Jacobi iteration for larger problems with more degrees of freedom.

We quantify the speedup achieved on each mesh using the SRJ method with heuristic and CJM relative to Jacobi iteration in Table \ref{tab:scaling}, which also shows the number of iterations required for each method to converge for each mesh. For a given geometry, the speedup achieved by both methods increases as the fidelity of the mesh increases. As an example, the lowest fidelity plate with hole mesh shows a nearly three times speedup when using SRJ, but this increases to nearly fourteen times when considering the finest plate with hole mesh. The CJM method outperforms the SRJ method for all meshes associated with the circle geometry. However, for the plate with hole and airfoil meshes, SRJ outperforms CJM and yields a higher speedup relative to Jacobi. For the finest airfoil mesh considered here, we expect a 25 times speedup relative to Jacobi iteration when using the SRJ method with heuristic (based on extrapolating the Jacobi iteration count). For this case, the CJM method takes twice as many iterations as SRJ (yielding an approximately 12 times speedup relative to Jacobi). 


\begin{table}[htbp!]
    \centering
    \begin{tabular}{|c|c|c|c|c c|c c|}
        \hline
        \textbf{Mesh} & \textbf{Fidelity} & \textbf{DOF Count} & \textbf{Jacobi} & \textbf{SRJ} & \textbf{Speedup} & \textbf{CJM} & \textbf{Speedup} \\ 
        \Xhline{1pt}
        \textit{Circle} & Low & 88 & 309 & 103 & 3.00 & 71 & 4.35  \\ 
        \hline
        & Medium & 362 & 1439 & 250 & 5.76 & 151 & 9.53  \\ 
        \hline
        & Fine & 1452 & 6095 & 571 & 10.67 & 299 & 20.38  \\ 
        \Xhline{1pt}
        \textit{Plate with Hole} & Low & 260 & 197 & 69 & 2.86 & 77  & 2.56  \\ 
        \hline
        & Medium & 1049 & 884 & 125 & 7.07 & 192 & 4.60  \\ 
        \hline
        & Fine & 4066 & 3593 & 260 & 13.82 & 514 & 6.99  \\ 
        \Xhline{1pt}
        \textit{Airfoil} & Low & 376 & 1290 & 199 & 6.48 & 554 & 2.33  \\ 
        \hline
        & Medium & 1630 & 5873 & 539 & 10.90 & 1172 & 5.01  \\ 
        \hline
        & Fine & 6106 & $\sim28470$ & 1133 & $\sim25.13$ & 2288 & $\sim12.44$  \\ 
        \hline
    \end{tabular}
    \caption{Comparison of the number of iterations required for Jacobi iteration, SRJ with the heuristic, and CJM to solve 2D Poisson on the three unstructured meshes (convergence is achieved when the $L_2$ residual norm reaches a value below $10^{-9}$). The CJM method shows convergence in the fewest iterations for the circle mesh (where the grid spacing is generally uniform). For the plate with hole and airfoil meshes, the SRJ method yields the best speedup. The airfoil mesh shows the most considerable reduction in number of iterations required for convergence, with approximately 25x speedup over Jacobi for the finest airfoil mesh.}
    \label{tab:scaling}
\end{table}

Our results illustrate that SRJ with the data driven heuristic is advantageous for more general unstructured PDE problems. The savings in computation relative to Jacobi increases as the fidelity of the mesh, or size of the linear system arising from discretization of the PDE is increased. The SRJ method with heuristic outperforms the CJM method particularly in cases where the mesh is highly nonuniform (comprised of very coarse and very fine grid spacing in different sections of the domain). This is particularly true for the airfoil mesh, where we observe very fine spacing near the airfoil, but much coarser spacing further away. In this case, it is difficult to determine an appropriate CJM scheme which captures all length scales well. The heuristic can still be used to determine SRJ schemes which give reasonable convergence on highly non-uniform meshes. 

In summary, the data based rule for selecting SRJ schemes provides a mechanism to apply SRJ effectively to a variety of matrices without the need to specifically tailor schemes to the given matrix or problem size. Despite being developed from convergence data from relatively small 1D Poisson matrices, the rule generalizes well to higher dimensional elliptic problems on both structured and unstructured meshes. For structured problems, the CJM method tends to give optimal convergence, although the user is required to determine this optimal scheme. For unstructured problems, particularly those which display a large degree of nonuniformity, the SRJ schemes are effective for convergence, whereas it may be difficult to determine an appropriate CJM scheme for efficient convergence.

\section{Conclusion}

The Scheduled Relaxation Jacobi method improves upon the convergence of the traditional Jacobi iteration by introducing relaxation parameters which effectively attenuate the solution error. In this paper, we have presented a family of schemes which can be used to achieve accelerated convergence for symmetric linear systems which would converge via Jacobi iteration (such as those arising from discretization of elliptic PDEs). Each scheme employs a different number of relaxation factors $M$ in a given cycle. It is desirable to find the scheme which gives optimal convergence when solving a linear system. While the best scheme for the optimal asymptotic convergence can be found analytically for certain problems, for general large scale problems, this may require experimentation with many different schemes to determine the one which yields the best convergence.

In this work, we have developed a data driven heuristic to determine which SRJ scheme in our family of schemes a user should employ during the linear solve process. The heuristic avoids the need to perform experimentation with many different SRJ schemes as it can automatically select a suitable scheme for the next cycle. Although the rule was developed using limited convergence data (specifically from applying SRJ schemes to linear systems arising from discretization of 1D Poisson), it outperforms an SRJ approach which uses a brute force level increasing rule and generalizes well to problems arising from discretization of PDEs in higher dimensions and on unstructured domains. The CJM schemes developed in previous work also exhibit good convergence, and in several cases outperform the SRJ schemes presented here. However, applying CJM can require user experimentation/analysis to determine a good scheme to utilize. Our heuristic removes the need to perform any analysis ahead of time. Additionally, the SRJ method with heuristic shows promising speedup capability for more general unstructured finite element problems which may be more representative of the large scale problems practitioners would like to solve.

The SRJ method provides a promising approach for solving large linear systems of equations. The data based heuristic provides an additional tool that allows practitioners to take advantage of the simplicity of SRJ and apply it very easily for solving a variety of problems without having to hand tailor schemes for each problem. One simply needs to derive SRJ schemes (given by Equation \eqref{eqn:srj-relaxation-parameters}) associated with different scheme levels provided in Appendix B, and allow the heuristic to decide an appropriate scheme level during the linear solve process. Furthermore, the method can provide good performance on the latest hardware architectures. Adsuara et al. have explored the performance of an SRJ implementation on GPUs \cite{Adsuara3}. Furthermore, implementation of SRJ schemes as smoothers within a larger multigrid framework can also provide additional convergence acceleration as shown by Yang and Mittal \cite{Yang17}. 

The SRJ algorithm may pave the way for a new class of high performance and parallel linear solvers for general PDE problems which utilize available simulation data to augment their capability and improve their ease of use.

\appendix

\clearpage
\section{Amplification polynomials as Chebyshev polynomials}

\begin{figure}[htbp!]
    \begin{subfigure}{.5\textwidth}
        \centering
        \includegraphics[width=.8\linewidth]{Figures/amplification_factors.png}  
        \caption{Amplification Polynomials $G_M(\lambda)$}
        \label{fig:sub-first}
    \end{subfigure}
    \begin{subfigure}{.5\textwidth}
        \centering
        \includegraphics[width=.8\linewidth]{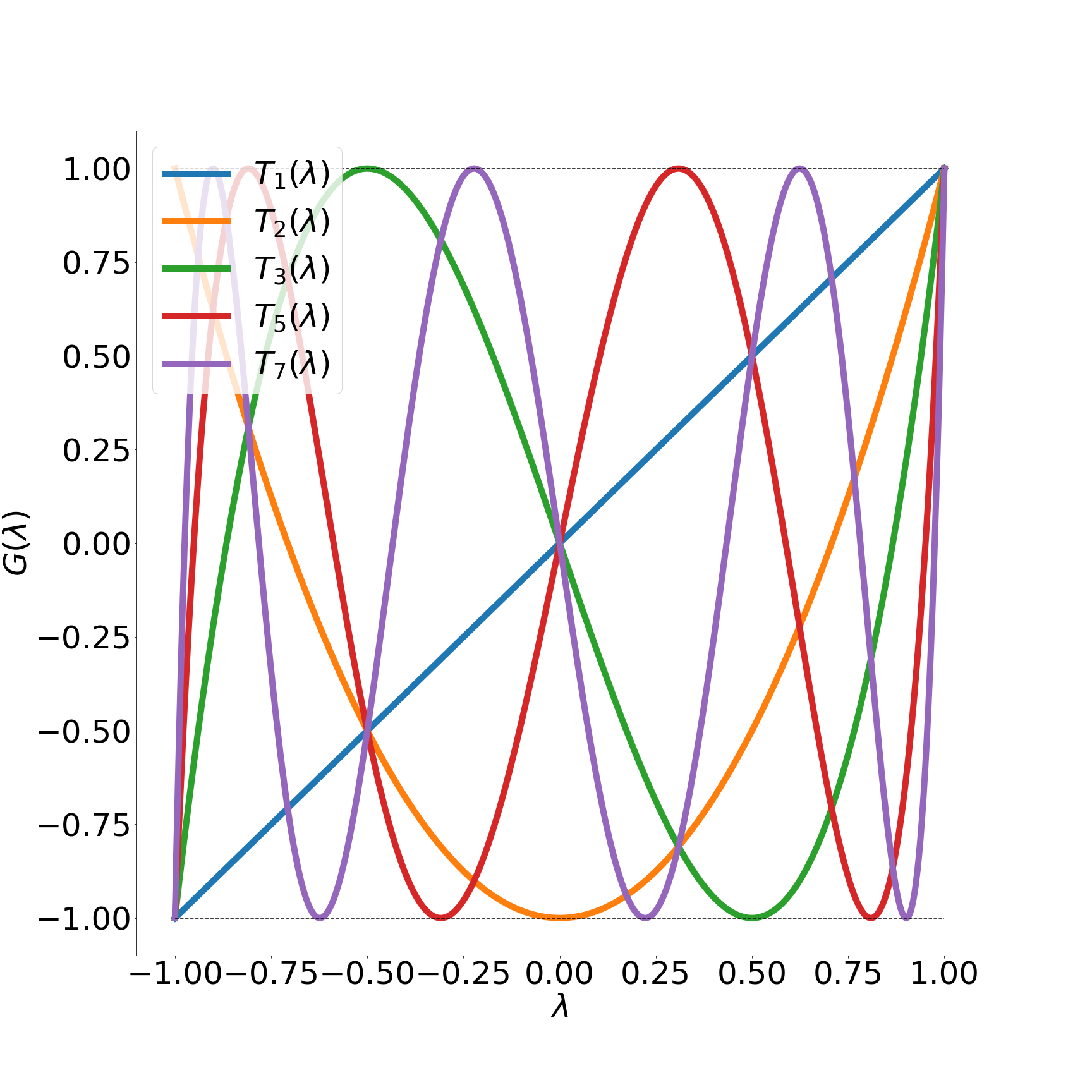}  
        \caption{Chebyshev polynomials $T_M(\lambda)$}
        \label{fig:sub-second}
    \end{subfigure}
     \caption{Amplification polynomials $G_{M}$ corresponding to SRJ schemes (left) and Chebyshev polynomials $T_{M}$ (right) for $M = 1,2,3,5,7$}
    \label{fig:amplification-chebyshev-polynomials}
\end{figure}

In this appendix, we derive the relationship between the amplification polynomials corresponding to our SRJ scheme and the Chebyshev polynomials, shown in Figure \ref{fig:amplification-chebyshev-polynomials} on the left and right respectively.

The original Chebyshev polynomials are bounded by $\pm 1$. However, the amplification polynomials we have derived are bounded by $\frac{1}{3}$, suggesting that a vertical scaling of $\frac{1}{3}$ is required to transform from the Chebyshev polynomials to our amplification polynomials. We can define an intermediate polynomial $\tilde{T}_{M} (\lambda) = \frac{1}{3} T_{M} (\lambda)$. In addition to a vertical scaling, a horizontal scaling is required to map between the two polynomials. In particular, we may define an affine transformation $f(\lambda)$ such that $G(\lambda) = \tilde{T}_M(f(\lambda))$. The affine transformation must satisfy two requirements which are derived below:
\begin{enumerate}
    \item By definition, $G(-1) = \tilde{T}_M(f(-1))$. Furthermore, it is true that  $G(-1) = \tilde{T}_M(-1)$. Therefore it follows that $G(-1) = \tilde{T}_M(f(-1)) = \tilde{T}_M(-1)$ so the affine transformation obeys $f(-1) = -1$.
    \item By definition, $G(1) =\tilde{T}_M(f(1))$. Define $\lambda^*$ as the argument which satisfies $T_{M}(\lambda^*) = 3$, so that $\tilde{T}_{M}(\lambda^*) = 1$. It is also true that $G_M(1) = 1$. Therefore, it follows that $G(1) =\tilde{T}_M(f(1)) = \tilde{T}_{M}(\lambda^*)$ so the affine transformation obeys $f(1) = \lambda^*$.
\end{enumerate}

The first condition follows from enforcing that the vertically scaled Chebyshev polynomial and the amplification polynomial have the same value at $\lambda = -1$, and the second condition ensures that the transformation will result in an amplification polynomial with value 1 at $\lambda = 1$. Assume that the affine transformation has the following form
\begin{equation}
    f(\lambda) = c_1 \lambda + c_0
\end{equation}
Substituting the two conditions $f(-1) = -1$ and $f(\lambda^*) = 1$ results in the following system of equations for the coefficients $c_0$ and $c_1$
\begin{equation}
    -c_1 + c_0 = -1 \\
\end{equation}
\begin{equation}
    c_1 + c_0 = \lambda^*
\end{equation}
Solving for the constants $c_0$ and $c_1$ of the affine transformation yields
\begin{equation}
    c_1 = \frac{\lambda^*+1}{2} \ , \   c_0 = \frac{\lambda^*-1}{2} 
\end{equation}
Therefore, the affine transformation is given by
\begin{equation}
    f(\lambda) = \frac{(\lambda^*+1) \lambda + (\lambda^*-1)}{2} 
    \label{eqn:affine}
\end{equation}
Lastly, the overall transformation which transforms the original Chebyshev polynomials $T_{M}(\lambda)$ to the amplification polynomials $G_{M}(\lambda)$ is
\begin{equation}
    G_{M}(\lambda) = \frac{T_{M}(f(\lambda))}{3}
\end{equation}
where $f(\lambda)$ is defined in Equation \eqref{eqn:affine}.

\clearpage
\section{Correspondence between Scheme Level and $M$}

In this appendix, we show the relationship between the SRJ scheme level and $M$ (given in Table \ref{tab:srj-scheme-level}). We chose SRJ schemes corresponding to specific $M$ to be selectable in our implementation, in order to prevent the available schemes from being too similar. The results in this paper can be reproduced by utilizing the specific SRJ schemes below.

\begin{table}[htbp!]
    \caption{Relationship between $M$ and scheme level used in this work}
    \centering
    \begin{tabular}{|c|c|}
    \hline
    Scheme Level & $M$ \\
    \hline
    0 &  1\\
    \hline
    1 &  2\\
    \hline
    2 & 3\\
    \hline
    3 & 5\\
    \hline
    4 & 7\\ 
    \hline
    5 & 10\\
    \hline
    6 & 14\\
    \hline
    7 & 19\\
    \hline
    8 & 26\\
    \hline
    9 & 35\\
    \hline
    10 & 47\\
    \hline
    11 & 63\\
    \hline
    12 & 84\\
    \hline
    13 & 111\\
    \hline
    14 & 147\\
    \hline
    15 & 194\\
    \hline
    16 & 256\\
    \hline
    17 & 338\\
    \hline
    18 & 446\\
    \hline
    19 & 589\\
    \hline
    20 & 778\\
    \hline
    21 & 1027\\
    \hline
    22 & 1356\\
    \hline
    23 & 1790 \\
    \hline
    24 & 2362 \\
    \hline
\end{tabular}
    \label{tab:srj-scheme-level}
\end{table}

\clearpage
\printbibliography

\end{document}